%% file: obstacle.tex
\input amstex
\input inputs

\def\rn{\Bbb{R}^n}
\def\obar{\overline{\Omega}}

\def\vbar{\overline{v}}
\def\wp{\omega^{+}}
\def\wm{\omega^{-}}

\def\l#1{{\cal L_{\rm#1}}}
\def\obar{\overline{\Omega}}

\def\vbar{\overline{v}}
\def\l#1{{\Cal L_{\text{\rm#1}}}}
\documentstyle{amsppt}
\hfuzz=9pt
\proofmodefalse
  
\topmatter
\title
THE OBSTACLE PROBLEM FOR FUNCTIONS OF LEAST GRADIENT
\endtitle
\author
William P. Ziemer\qquad
Kevin Zumbrun
\endauthor
\dedicatory 
Dedicated to Alois Kufner on the occasion of his
sixty-fifth birthday, with the authors' best wishes.
\enddedicatory  
\address
Department of Mathematics, Indiana University, Bloomington, IN  47405--5701
\endaddress
\thanks Research of first author supported in part while visiting
Charles University, June, 1997\endthanks
\thanks Research of second author supported in part by a grant from
the National Science Foundation\endthanks 
 
 
\abstract
For a given domain $\Omega \subset \Bbb{R}^n$, we consider the variational
problem of minimizing the $L^1$-norm of the gradient on $\Omega$ of a function
$u$ with prescribed continuous boundary values and satisfying a continuous
lower obstacle condition $u\ge \Psi$ inside $\Omega$.
Under the assumption of strictly positive mean curvature of
the boundary $\partial\Omega$, we show existence of a continuous
solution, with H\"older exponent half of that of data and obstacle.

This generalizes previous results obtained for the unconstrained and 
double-obstacle problems.  The main new feature in the present analysis 
is the need to extend various maximum principles from the case of two
area-minimizing sets to the case of one sub- and one superminimizing
set.  This we accomplish subject to a weak regularity assumption
on one of the sets, sufficient to carry out the analysis. Interesting
open questions include the uniqueness of solutions and a complete
analysis of the regularity properties of area superminimizing 
sets.  We provide some preliminary results in the latter direction, namely
a new monotonicity principle for superminimizing sets, and the
existence of ``foamy'' superminimizers in two dimensions.
\endabstract
\endtopmatter
       
\document
\chapno=1
\subsectionnumber=0
\Increase{\sectionnumber}
\head
\the\chapno. Introduction
\endhead
\thmno=0
\tagno=0 

A rather complete and extensive literature is now in place concerning
existence and regularity of solutions to a wide range of variational
problems for which the following is prototypical: 
$$
\inf\left\{\int_{\Omega }\abs{\nabla u}^{p}:u\in W^{1,p}(\Omega ),u-g\in
W^{1,p}_{0}(\Omega )\right\}.\eqnlbl{pdirichlet}
$$
Here, $\Omega \subset\rn$ is a bounded, open set,
$1<p<\infty$ and $g\in W^{1,p}(\Omega )\cap C^{0}(\overline{\Omega
})$. The Euler-Lagrange equation for \eqnref{pdirichlet} is the
$p$-Laplacian $\div(\abs{\nabla u}^{p-2}\nabla u)=0$. The interested
reader can consult recent books on this subject and the references
therein, \cite{AH}, \cite{HKM}, and \cite{MaZ}. The theory related to the case
corresponding to $p=1$ is far less complete. In spite of the fact
that there is a vast literature relating to the least area functional,
$$
\inf_{u}\left\{\int_{\Omega }\sqrt{1+\abs{\nabla u}^{2}}\right\},
$$
there are many open questions concerning other functionals with
linear growth in $\abs{\nabla u}$. Investigations concerning such
questions were considered in [SWZ], [SZ], ???. In particular, the Dirichlet problem was
investigated; that is, for a bounded Lipschitz domain
$\Omega \subset \Bbb{R}^n$, and for $g:\partial \Omega \to \Bbb{R}^1$
continuous, the questions of existence and regularity of solutions to 
$$
\inf \left\{ \|\nabla u\|(\Omega): u\in BV(\Omega),\ u=g\ on\
\partial\Omega\right\} \eqnlbl{1.2}
$$
were examined. Here $\|\nabla u\|(\Omega)$ denotes the total variation of the
vector-valued measure $\nabla u$ evaluated on $\Omega$.  It was shown
that a solution $u\in BV(\Omega)\cap C^0(\Omega)$ exists provided
that $\partial\Omega$ satisfies two conditions, namely, that
$\partial\Omega$ has non-negative curvature (in a weak sense) and that
$\partial\Omega$ is not locally area-minimizing. See Section 2 below for
notation and definitions. 

In this paper we consider the obstacle problem 
$$
\inf \left\{\norm{\nabla u}(\Omega ) : u\in C^0(\overline{\Omega}),\
u\geq\psi \ on\; \Omega, u=g\geq\psi \ on\
\partial\Omega  \right\}\eqnlbl{1.3}
$$
where $\Omega \subset\rn$ is a bounded Lipschitz domian,
$g\colon\partial \Omega \to\Bbb{R}^1$ is continuous and $\psi $ is a
continuous function on $\overline\Omega$. The analogous obstacle
problem for \eqnref{pdirichlet} was investigated by several authors
and is now well understood, cf. \cite{CL}, \cite{Li}, \cite{MiZ}, \cite{MuZ}. One of the
difficulties encountered in the analysis of both \eqnref{pdirichlet}
and \eqnref{1.3} is the fact that the compactness in $L^1(\Omega)$ of
a sequence whose $BV$-norms are bounded does not ensure, a {\it
priori}, continuity of the limiting function or that it will assume
the boundary values $g$, thus making the question of existence
problematic. In this paper as well as in \cite{SWZ}, we rely heavily on
the discovery made in \cite{BDG} that the
superlevel sets of a function of least gradient are area-minimizing.
This fact, along with the co-area formula (see 
\eqnref{2.10} below), suggests that the existence of a function of least
gradient subject to an obstacle constraint can be established by actually
constructing each of its superlevel sets in such a way that it reflects
both the appropriate boundary condition and the obstacle condition. The main thrust of this paper
is to show that this is possible. Thus we show that there exists a continuous
solution to \eqnref{1.3} and we also show it inherits essentially the same
regularity as the boundary data and obstacle.

As in \cite{SWZ}, both existence and regularity are developed by extensive
use of $BV$ theory and sets of finite perimeter as well as certain
maximum principles. One of the main contributions of this paper is a
new maximum principle that involves a super area-minimizing set
and an area-minimizing set, Theorem \thmref{maxprin}. The similar
result involving two area-minimizing sets, due independently to \cite{Mo}
and \cite{S2}, played a crucial role in \cite{SWZ}.

Our extended maximum principle requires a weak regularity property
on one of the sets, that the set be contained in
the (topological) closure of its interior.
This is clearly satisfied in the contexts that we apply it, for which one
of the sets is always area-minimizing.  But, an interesting open question
is whether or not this technical assumption can be dropped.  

This issue leads us to consider a question of interest in its own right:
``What is the regularity of a (sub)superminimizing set?''  We conclude by
presenting some separate,  preliminary results on this subject, including
a new monotonicity principle for (sub)superminimizing sets, and the
existence of unusual, ``foamy'' (sub) superminimizers in two dimensions.
It is our hope that these results will stimulate further investigation 
into the topic of regularity.

\chapno=2
\subsectionnumber=0
\Increase{\sectionnumber}
\head
\the\chapno. Preliminaries
\endhead
\thmno=0
\tagno=0 

The Lebesgue measure of a set $E\subset \Bbb{R}^n$ will be denoted $|E|$ and
$H^\alpha(E), \alpha>0$, will denote $\alpha$-dimensional Hausdorff measure of
$E$.  Throughout the paper, we almost exclusively employ $H^{n-1}$.  The
Euclidean distance between two points $x,y\in\Bbb{R}^n$ will be denoted by
$|x-y|$. The open ball of radius $r$ centered at $x$ is denoted by
$B(x,r)$ and $\bar{B}(x,r)$ denotes its closure.

If $\Omega\subset \Bbb{R}^n$ is an open set, the class of function $u\in
L^1(\Omega)$ whose partial derivatives in the sense of distribution are
measures
with finite total variation in $\Omega$ is denoted by $BV(\Omega)$ and is
called the space of functions of bounded variation on $\Omega$.  The space
$BV(\Omega)$ is endowed with the norm
 
$$
\|u\|_{BV(\Omega)} = \|u\|_{1;\Omega} + \|\nabla u\|(\Omega) \eqnlbl{2.1}
$$
 
where $\|u\|_{1;\Omega}$ denotes the $L^1$-norm of $u$ on $\Omega$ and
where 
$\|\nabla u\|$ is the total variation of the vector-valued measure $\nabla
u$. 

The following compactness result for $BV(\Omega)$ will be needed later, cf. \cite{Gi}
or \cite{Z}.

\proclaim{Theorem 2.1}  If $\Omega\subset \Bbb{R}^n$ is a bounded Lipschitz
domain, then
 
$$
BF(\Omega)\cap \left\{ u:\|u\|_{BV(\Omega)}\leq 1\right\} 
$$
 
is compact in $L^1(\Omega)$.  Moreover, if $u_i\to u$ in $L^1(\Omega)$ and
$U\subset \Omega$ is open, then
 
$$
\liminf_{x\to\infty} \|\nabla u_i\|(U) \geq \|\nabla u\|(U)
$$
\endproclaim

A Borel set $E\subset \Bbb{R}^n$ is said to have {\it finite perimeter in}
$\Omega$ provided the characteristic function of $E, \chi_E$, is a function of
bounded variation in $\Omega$.  Thus, the partial derivatives of $\chi_E$ are
Radon measures on $\Omega$ and the perimeter of $E$ in $\Omega$ is defined as
 
$$
P(E,\Omega) = \|\nabla\chi_E\|(\Omega). \eqnlbl{2.2}
$$

A set $E$ is said to be of {\it locally finite perimeter} if $P(E,\Omega) <
\infty$ for every bounded open set $\Omega\subset \Bbb{R}^n$.

One of the fundamental results in the theory of sets of finite perimeter is
that
they possess a measure-theoretic exterior normal which is suitably general to
ensure the validity of the Gauss-Green theorem.  A unit vector $\nu$ is defined
as the measure-theoretic exterior normal to $E$ at $x$ provided
$$
\lim_{r\to 0} r^{-n} |B(x,r) \cap \left\{ y:(y-x)\cdot\nu < 0,\ y\not\in
E\right\} |= 0
$$
and
$$
\lim_{r\to 0} r^{-n} | B(x,r) \cap \left\{ y:(y-x)\cdot\nu > 0,\ y\not\in
E\right\} |=0. \eqnlbl{2.3}
$$
The
measure-theoretic normal of $E$ at $x$ will be denoted by $\nu(x,E)$ and we
define
$$
\partial_*E = \left\{ x:\nu(x,E)\ exists \right\}.
\eqnlbl{2.4}
$$
The Gauss-Green theorem in this context 
states that if $E$ is a set of locally finite perimeter and 
$V\colon R^{n}\to R^{n}$ is a Lipschitz vector field, then 
$$
\int_{E}\div
V(x)\,dx=\int_{\partial_{*}E}V(x)\cdot\nu(x,E)\,dH^{n-1}(x),\eqnlbl{gg} 
$$
cf. \cite{Fe2, \S 4.5.6}.
Clearly, $\partial_*E\subset \partial E$, where $\partial E$ denotes
the topological boundary of $E$.  Also, the {\it topological interior
of} $E$ is denoted by $E^i = (\Bbb{R}^n \setminus \partial E) \cap
E$, the {\it topological exterior} by $E^e = (\Bbb{R}^n \setminus
\partial E) \cap (\Bbb{R}^n \setminus E)$ and $E^c$ to denote the
complement $\Bbb{R}^n\setminus E$.  The notation $E\subset\subset F$
means that the closure of $E,\;\overline E$, is a compact subset of $F^i$.  

For measurable sets $E$, the {\it measure-theoretic interior},
$E_{m}^{i}$, is the set of all points at which the metric density of
$E$ is 1 and the {\it measure-theoretic exterior}, $E_{m}^{e}$, is
all points at which the metric density is $0$. The measure
theoretic-boundary, $\partial _{m}E:=\Bbb{R}^n\setminus(E_{m}^{i}\cup
E_{m}^{e})$. Clearly, $\partial_*E \subset
\partial_{m}E\subset \partial E$.  Moreover, it is well known that
$$
E\ \text{{\it is of finite perimeter if and only if}}\ 
 H^{n-1}(\partial_{m}E)<\infty \eqnlbl{2.6}
$$
and that
$$
P(E,\Omega) = H^{n-1} (\Omega \cap \partial_{m}E) = H^{n-1} (\Omega \cap
\partial_*E)\ \text{{\it whenever}}\ P(E,\Omega)<\infty \eqnlbl{2.7}
$$
cf. \cite{Fe2 \S4.5}.  From this it easily follows that 
$$
P(E\cup F,\Omega) + P(E\cap F,\Omega) \leq P(E,\Omega) + P(F,\Omega),
\eqnlbl{2.8}
$$
thus implying that sets of finite perimeter are closed under finite unions and
intersections.

The definition implies that sets of finite perimeter are defined only up to
sets of measure 0.  In other words, each set determines an equivalence class of
sets of finite perimeter.  In order to avoid this ambiguity, we will
employ 
$$
\tilde E:=E\cup E_{1}\setminus E_{0} \eqnlbl{defineE}
$$
as the distinguished representative for $E$. Here, 
$$
E_{i}:=\{x:\abs{E\cap B(x,r)}/\abs{B(x,r)}=i\;\text{for all small
$r>0$}\}, \;i=0,1.
$$
Thus, with this convention, it easy to see that 
$$
\overline{\partial_*E} = \partial E. \eqnlbl{2.10}
$$
This convention will apply, in particular, to all competitors of the
variational problems \eqnref{2.22} and \eqnref{2.23} below as well as to the
sets defined by \eqnref{2.18}.

Of particular importance to us are sets of finite perimeter whose
boundaries are area-minimizing.  If $E$ is a set of locally finite
perimeter and $U$ a bounded, open set, then $E$ is said to be {\it
area-minimizing in $U$} if $P(E,U)\leq P(F,U)$ whenever $E\Delta
F\subset\subset U$. Also, $E$ is said to be {\it super area-miniming
in $U$} ({\it sub area-minimizing in $U$}) if $P(E,U)\leq P(E\cup F,U)$
($P(E,U)\leq P(E\cap F,U)$) whenever $E\Delta F\subset\subset U$.

A tool that will play a significant role in this paper is the co-area
formula.  It states that if $u\in BV(\Omega)$, then
$$
\|\nabla u\|(\Omega) = \int^\infty_{-\infty} P(E_t,\Omega)dt \eqnlbl{coarea}
$$
where $E_t = \{u\geq t\}$.  In case $u$ is Lipschitz, we have
$$
\int_\Omega |\nabla u|dx = \int^\infty_{-\infty} H^{n-1} \left( u^{-1}(t) \cap
\Omega\right) dt.
$$
Conversely, if $u$ is integrable on $\Omega$ then
$$
\int^\infty_{-\infty}P(E_t,\Omega)dt < \infty\ \text{{implies}}\ u\in
BV(\Omega),
\eqnlbl{reversecoa}
$$
cf. \cite{Fe1}, \cite{FR}.

Another fundamental result is the isoperimetric inequality for sets
of finite perimeter. It states that there is a constant $C=C(n)$ such
that 
$$
P(E)^{n/(n-1)}\leq C\abs{E}\eqnlbl{iso}
$$
whenever $E\subset\rn$ is a set of finite perimeter. Furthermore,
equality holds if and only if $E$ is a ball.

The regularity of $\partial E$ plays a crucial role in our development.  In
particular, we will employ the notion of tangent cone.  Suppose $E$ is
area-minimizing in $U$ and for convenience of notation, suppose $0\in U\cap
\partial E$.  For each $r>0$, let, $E_r = \Bbb{R}^n \cap \{x:rx\in E\}$.  It is
known (cf. \cite{S1, \S35}) that for each sequence $\{r_i\}\to 0$,
there exists a subsequence (denoted by the full sequence) such that
$\chi_{E_i}$ converges in $L^1_{loc}(\Bbb{R}^n)$ to $\chi_C$, where $C$ is a
set of locally finite perimeter.  In fact, $C$ is area-minimizing and
is called the tangent cone to $E$ at 0.  Although it is not immediate, $C$ is a
cone and therefore the union of half-lines issuing from 0.  It follows from [S1,
\S37.6] that if $\overline{C}$ is contained in $\overline{H}$ where $H$ is any
half-space in $\Bbb{R}^n$ with $0\in\partial H$, then $\partial H$ is regular
at).  That is there exists $r>0$ such that 
 
$$
B(0,r)\cap \partial E\ \text{ is a real analytic
hypersurface.}\eqnlbl{2.14}
$$
 
Furthermore, $\partial E$ is regular at all points of $\partial_*E$ and
 
$$
H^\alpha\left( (\partial E\setminus \partial_*E)\cap U\right) = 0\ \text{{\it for all}}\
\alpha > n-8, \eqnlbl{2.15}
$$
 
cf. \cite{Gi, Theorem 11.8}.

The boundary data $g$ admits a continuous extension $G\in
BV(\Bbb{R}^n\setminus \overline{\Omega}) \cap C^0(\Bbb{R}^n\setminus
\Omega)$, \cite{Gi, Theorem 2.16}.  In fact, $G\in C^\infty(\Bbb{R}^n
\setminus \overline{\Omega})$, but we only need that $G$ is
continuous on the complement of $\Omega$.  Clearly, we can require
that the support of $G$ is contained in $B(0,R)$ where $R$ is chosen
so that $\Omega\subset\subset B(0,R)$.  We have
$$
G\in BV(\Bbb{R}^n \setminus  \overline{\Omega}) \cap C^0(\Bbb{R}^n \setminus  \Omega)\ \text{{\it
with}}\ G = g\ \text{{\it on}}\ \partial\Omega. \eqnlbl{2.17}
$$

We now introduce sets that will ensure that our constructed solution
satisfies the required Dirichlet condition $u = g$ on
$\partial\Omega$ and the obstacle condition $u\geq \psi $ in
$\Omega$.  Thus, for each $t\in [a,b]$, let
$$
\Cal{L}_t = (\Bbb{R}^n \setminus  \Omega) \cap \{x:G(x) \geq t\},\eqnlbl{2.18}
$$
 
$$
L_{t} = \text{closure}(\left\{ x:x\in\Omega, \psi (x)>t \right\}). \eqnlbl{2.19}
$$
 
Note that the co-area formula \eqnref{coarea} and the fact that $G\in
BV(\Bbb{R}^n\setminus\overline{\Omega})$ imply that
$P(\Cal{L}_t,\Bbb{R}^n\setminus \overline{\Omega})<\infty$ for almost all $t$.  For all
such $t$, we remind the reader that we employ our convention \eqnref{defineE} in
defining $\Cal{L}_t$.

We let $[a,b]$ denote the smallest interval containing $g(\partial
\Omega )\cup\psi (\overline\Omega )$ and define
$$
T := [a,b]\cap \left\{
t:P(\Cal{L}_t,\Bbb{R}^n\setminus \overline{\Omega})<\infty.\right.\eqnlbl{2.20}
$$
Thus, by \eqnref{2.7} and the fact that $H^{n-1}(\partial \Omega)<\infty$, we obtain
$$
H^{n-1}(\partial_{m}\Cal{L}_t) = P(\Cal{L}_t,\Bbb{R}^n \setminus  \overline{\Omega}) +
H^{n-1} \left[(\partial_{m}\Cal{L}_t)\cap (\partial \Omega)\right] < \infty.
\eqnlbl{2.21}
$$
 
For each $t\in T$, the variational problems
 
$$
\min\left\{ P(E,\Bbb{R}^n): E \setminus  \overline{\Omega} = \Cal{L}_t \setminus 
\overline{\Omega}, \;E\supset L_t\right\}, \eqnlbl{2.22}
$$
 
$$
\max \left\{ |E|:E\;\;\text{is a solution of \eqnref{2.22}}\right\}
\eqnlbl{2.23}
$$
will play a central role in our development.  In light of Theorem 2.1, a
solution to both problems can be obtained from the direct method.  \eqnref{2.21}
is also used to obtain existence for \eqnref{2.22}.  We will denote by $E_t$ the
solution to \eqnref{2.23}.  In this regard, note that our convention
\eqnref{defineE} ensures that $E_t \setminus  \overline{\Omega} = \Cal{L}_t \setminus 
\overline{\Omega}$; furthermore, because of our convention, $\Cal{L}_t$ need
not be a closed set. Also, observe that that $E_{t}$ is super
area-minimizing in $\Omega $.

\chapno=3
\subsectionnumber=0
\Increase{\sectionnumber}
\head
\the\chapno. A Maximum Principle
\endhead
\thmno=0
\tagno=0

First, we begin with a result which is a direct consequence of a maximum
principle for area-minimizing hypersurfaces established independently
in \cite{Mo} and \cite{S2}.

\proclaim{\thmlbl{Theorem 2.2} Theorem}  Let $E_1\subset E_2$ 
and suppose both $E_{1}$ and $E_{2}$ are area-minimizing in an open
set $U\subset \Bbb{R}^n$.  Further, suppose $x\in (\partial E_1)\cap
(\partial E_2) \cap U$.  Then $\partial E_1$ and $\partial E_2$ agree
in some neighborhood of $x$.
\endproclaim

\proclaim{\thmlbl{Lemma A} Lemma} For arbitrary measurable sets
$A,B\subset\Bbb{R}^n$, it holds that 
$$\align
H^{n-1}(\partial_{m}(A\cup B))
&\leq H^{n-1}(\partial_{m}A\cap (B_m^i)^c) +H^{n-1}(\partial_{m}B\cap (\bar A_m)^c)\\
H^{n-1}(\partial_{m}(A\cap B))
&\leq H^{n-1}(\partial_{m}A\cap B^i_m)+H^{n-1}(\partial_{m}B\cap \bar A_m).
\endalign
$$
\endproclaim
\demo{Proof}
It follows immediately from definitions that 
$$
(\partial_{m}A \cap (B^i_m)^c) \cup
(\partial_{m}B \cap (A^i_m)^c)  
=
(\partial_{m}A \cap (B^i_m)^c) \cup
(\partial_{m}B \cap (\bar A_m)^c),
$$
which yields the first inequality. 
The result for intersections then follows from $\partial_{m}A= \partial_{m}A'$
and $A\cap B= (A^c\cup B^c)^c$.
\enddemo
\qqed
\proclaim{\thmlbl{maxprin} Theorem} Let $E$ be sub area-minimizing  and
$F$ super area-minimizing relative to
an open set $U$, with $E\subset F$ and $\partial E\cap\partial
F\subset\subset U$. Further, suppose that $\overline{E}\cap
U=\overline{E^{i}}\cap U$.
Then, relative to $U$, 
either $\partial E \cap \partial F= \emptyset$  
or else $\partial E=\partial F$ in a neighborhood of $\partial E \cap
\partial F$. 
\endproclaim
\demo{Proof}
Suppose $\partial  E\cap \partial F\not= \emptyset$.
The set $\partial E\cap \partial F$ is 
contained in open neighborhood $V\subset \subset U$ and thus,
for sufficiently small $|w|,\;w\in\rn$, we have
$$
(E+w) \setminus F \subset V+ w \subset \subset U. 
\eqnlbl{1}
$$
Choose $x_0 \in \partial E\cap \partial F$. Since $E=\overline{E^{i}}$,
there exists $w\in\rn$ with  $|w|$ arbitrarily small
such that $x_0 - w\in E^i$, or equivalently
$$
x_0 \in  (E+w)^i.
\eqnlbl{2}
$$
Denote the translated set
$E+w$ by $E_w$.  By shrinking $U$ if necessary, we can arrange
that $E_w$ is sub area-minimizing in $U$.

Now we will show that $F$ is area-minimizing in
the open set $ U\cap E_w^i$.  For, suppose to
the contrary that there were a set $G$ with 
$$
G\Delta F \subset \subset  U\cap E_w^i
\eqnlbl{3}
$$
and
$$
P(G,  U\cap E_w^i) < P(F,  U\cap E_w^i).
\eqnlbl{4}
$$
By \eqnref{3}, $G\cap U=F\cap U$ near $\partial E_w$,
while, by \eqnref{4},
$$H^{n-1}(\partial_m G \cap E_w^i\cap U)<
H^{n-1}(\partial_m F \cap E_w^i\cap U)).
\eqnlbl{5}
$$ 
Since $F$ and $G$ agree on $(E_{w})^{i}_{m}\setminus
E_{w}^{i}\subset\partial_{m} E_{w}$, it follows that 
$$
H^{n-1}(\partial_m G \cap (E_{w})^{i}_{m}\cap U)<
H^{n-1}(\partial_m F \cap (E_{w})^{i}_{m}\cap U)).
\eqnlbl{5a}
$$ 
On the other hand, super area-minimality of $F$ in $U$
implies that $P(F\cup E_w,U)\geq P(F,U)$.  With
Lemma \thmref{Lemma A}, this gives
$$
\eqalign{
H^{n-1}(\partial_m F\cap {((E_w)^i_m)^c}\cap U)& +
H^{n-1}(\partial_mE_w\cap (\bar F_m)^c\cap U) \cr
&\geq
H^{n-1}(\partial_m (F\cup E_w)\cap U) \cr
&\geq 
H^{n-1}(\partial_m F\cap U) \cr
&= H^{n-1}(\partial F\cap {((E_w)^i_m)^c}\cap U)\cr
&\qquad +
H^{n-1}(\partial F\cap (E_w)^i_m\cap U), 
}
$$
and thus
$$
H^{n-1}(\partial E_w\cap (\bar F_m)^c\cap U) 
\geq
H^{n-1}(\partial F\cap (E_w)^i_m\cap U).
\eqnlbl{6}
$$
Therefore,
$$
\eqalign{
H^{n-1}(\partial(G\cap E_w) \cap U) 
&\leq
H^{n-1}(\partial E_w \cap \bar G_m \cap U) + 
H^{n-1}(\partial G \cap (E_w)^i_m \cap U) \cr
&=
H^{n-1}(\partial E_w \cap \bar F_m\cap U) + 
H^{n-1}(\partial G \cap (E_w)^i_m\cap U) \cr
&< H^{n-1}(\partial E_w \cap \bar F_m\cap U) +
H^{n-1}(\partial F \cap (E_w)^i_m\cap U) \cr
&\leq
H^{n-1}(\partial E_w \cap \bar F_m\cap U)+
H^{n-1}(\partial E_w \cap (\bar F_m)^c\cap U) \cr
&=
H^{n-1}(\partial E_w\cap U), 
}
\eqnlbl{7}
$$
where the first inequality follows by Lemma \thmref{Lemma A}, the second
by substituting $F$ for $G$ in the vicinity of $\partial E_w$,
the third by \eqnref{5}, the fourth by \eqnref{6}, 
and the last by set decomposition.
In other words, $P(G\cap E_w,U)< P(E_w,U)$.
But, at the same time, 
$$(G\cap E_w)\Delta E_w = E_w \setminus G
\subset (E_w\setminus F) \cup (G \Delta F)
$$
is compactly supported in $U$, by \eqnref{1} and \eqnref{3}.
contradicting the sub area-minimality of $E_w$ in $U$.
By contradiction, we have that $F$ is area-minimizing in
$E_w^i\cap U$, as claimed. 

By basic regularity results, we thus have also that $\overline(F)=
\overline{F^i}$ in a neighborhood of $x_0$.
By a symmetric argument, it follows that $E$ is
area-minimizing near $x_{0}$  as well, and therefore we can 
appeal to Theorem \thmref{Theorem 2.2} to obtain our conclusion.
\qqed
\enddemo

We do not know whether the hypothesis $E\cap U=\overline{E^{i}}\cap
U$ in the previous result is necessary. However, in the case where
$E$ is area-minimizing in $U$, the regularity results \eqnref{2.15}
show that the hypothesis is satisfied and this is sufficient for the
purposes of this paper. The following result is what we need and it
now follows immediately from Theorem \thmref{maxprin}. 

\proclaim{\thmlbl{maxprinsuper} Corollary} Let $E$ be area-minimizing  and
$F$ super area-minimizing relative to
an open set $U$, with $E\subset F$ and $\partial E\cap\partial
F\subset\subset U$. Then, relative to $U$,  
either $\partial E \cap \partial F= \emptyset$  
or else $\partial E=\partial F$ in a neighborhood of $\partial E \cap
\partial F$. 
\endproclaim

\chapno=4
\thmno=0
\tagno=0
\subsectionnumber=0
\Increase{\sectionnumber}
\head
\the\chapno. Construction of the solutuion
\endhead
In this section we will construct a solution $u$ of \eqnref{1.3} by
using $E_t \cap \overline{\Omega}$ to define the set $\{u \geq t\}$
up to a set of measure zero for almost all $t$.  This construction
will be possible for bounded Lipschitz domains $\Omega$ whose
boundaries satisfy the following two conditions.  \medskip

(i)  For every $x\in \partial\Omega$ there exists $\varepsilon _0>0$ such that for
every set of finite perimeter $A\subset\subset B(x,\varepsilon _0)$
 
$$
P(\Omega,\Bbb{R}^n) \leq P(\Omega \cup A, \Bbb{R}^n). \eqnlbl{3.1}
$$

(ii)  For every $x\in \partial\Omega$, and every $\varepsilon  \geq 0$ there exists
a set of finite perimeter $A\subset\subset B(x,\varepsilon )$ such that
 
$$
P \left(\Omega,B(x,\varepsilon )\right) > P\left( \Omega \setminus  A, B(x,\varepsilon )\right).
\eqnlbl{3.2}
$$
Clearly, we may assume that $x\in A$.

The first condition states that $\partial\Omega$ has non-negative mean
curvature (in the weak sense) while the second states that $\Omega$ is
not locally area-minimizing with respect to interior variations.  Also, it can
be easily verified that if $\partial\Omega$ is smooth, then both conditions
together are equivalent to the condition that the mean curvature of
$\partial\Omega$ is positive on a dense set of $\partial\Omega$.

Since $\Omega$ is a Lipschitz domain, for each $x_0 \in \partial \Omega,\
\partial\Omega$ can be represented as the graph of a nonnegative Lipschitz
function $h$ defined on some ball $B(x'_0,r)\subset \Bbb{R}^{n-1}$ where
$x'_0\in \Bbb{R}^{n-1}$.  That is, $\left\{ \left( x',h(x')\right):x'\in
B'(x'_0,r)\right\}$ $\subset \partial\Omega$.  Throughout we will use the
notation $B'(x'_0,r)$ and $X'$ to denote elements in $\Bbb{R}^{n-1}$ and thus
they will be distinguished them from their $n$-dimensional counterparts
$B(x_0,r)$ and $x$.

We assume our configuration is oriented in such a way that $$\left\{
(x',x''):0<x''<h(x')\right\} \subset \Omega.$$ With $S = \left\{
\left( x',h(x')\right):x'\in B'(x'_0,r)\right\}$ we have that
$$
H^{n-1}(S) = \int_{B'(x'_0,r)} \sqrt{1+|\nabla h|^2}\, dH^{n-1}(x'). 
$$
These facts lead immediately to the following result.
\proclaim{\thmlbl{lemma3.1} Lemma}  If $\Omega$ is a Lipschitz domain with non-negative mean
curvature in the sense of \eqnref{3.1}, then the function $h$, whose graph
represents $\partial\Omega$ locally, is a weak supersolution of the minimizing
surface equation.  That is, for $r$ sufficiently small,
$$
\int_{B'(x'_0,r)} \frac{\nabla h\cdot \nabla\phi}{\sqrt{1+|\nabla h|^2}}\,
dx'\geq 0
$$
whenever $\phi\in C^{1,1}_0 \left( B'(x'_0,r)\right),\phi>0$.
\endproclaim
We will also need the following result from \cite{SWZ, Lemma 4.2} whose proof is
an easy consequence of the weak Harnack inequality.
\proclaim{\thmlbl{Lemma 3.2} Lemma}
Suppose $W$ is an open subset of $\Bbb{R}^{n-1}$.  If
$v_1,v_2\in C^{0,1}(W)$ are respectively weak super and subsolutions of the
minimal surface equation in $W$ and if $v_1(x'_0) = v_2(x'_0)$ for some $x'_0
\in W$ while $v_1(x')\geq v_2(x')$ for all $x'\in W$, then
$$
v_1(x') = v_2(x')
$$
for all $x'$ in some closed ball contained in $W$ centered at $x'$.
\endproclaim

An important step in our development is the following lemma.

\proclaim{\thmlbl{scontain} Lemma}  For almost all $t\in [a,b],\; \partial E_t\cap \partial
\Omega\subset g^{-1}(t)$.
\endproclaim

\demo{Proof}  First note that if $t>\max_{x\in \partial\Omega} g(x)$, then
$\partial E_t\cap \partial \Omega = \emptyset$.  So we may assume that $t\in T$
and $t\leq \max_{x\in \partial\Omega}g(x)$.  The proof will proceed by
contradiction and we first show that $\partial E_t$ is locally area minimizing
in a neighborhood of each point $x_0 \in \left( \partial E_t \cap \partial
\Omega \setminus  g^{-1}(t)\right)$, i.e., we claim that there exists $\varepsilon  >0$, such
that for every set $F$ with the property that $F\Delta E_t \subset\subset
B(x_0,\varepsilon )$, we have
$$
P\left( E_t,B(x_0,\varepsilon )\right) \leq P\left( F,B(x_0,\varepsilon )\right)
\eqnlbl{3.3}
$$
or equivalently, $P(E_t,\Bbb{R}^n) \leq P(F,\Bbb{R}^n)$.

By our assumption, either $g(x_0) < t$ or $g(x_0)>t$.  First consider
the case $g(x_0)<t$.  Since $G(x_0) = g(x_0) < t$ and $G$ is
continuous on $\Bbb{R}^n\setminus \Omega$, there exists $\varepsilon >0$, such
that $B(x_0,\varepsilon )\cap \Cal{L}_t = \emptyset$.  Also, $\psi $ is
continuous on $\overline{\Omega}$ and $\psi (x_0) \leq g(x_0)<t$, so we
may take $\varepsilon$ small enough such that $L_t\cap
B(x_0,\varepsilon ) = \emptyset$.  We will assume that $\varepsilon
<\varepsilon _0$, where $\varepsilon_0$ appears in condition
\eqnref{3.1}.  We proceed by taking a variation $F$ that satisfies
$F\Delta E_t \subset\subset B(x_0,\varepsilon )$.  Because
of \eqnref{3.1} and 
$$
P(E\cup F,\Omega) + P(E\cap F,\Omega) \leq
P(E,\Omega) + P(F,\Omega)
\eqnlbl{2.8a}, note that 
$$
for every $A\subset\subset B(x_0,\varepsilon )$,
$$
\aligned
P(A\cup \Omega,\Bbb{R}^n) + P(A\cap \Omega,\Bbb{R}^n) & \leq P(A,\Bbb{R}^n) +
P(\Omega,\Bbb{R}^n) \\
& \leq P(A,\Bbb{R}^n) + P(A\cup \Omega,\Bbb{R}^n) 
\endaligned
$$
 
Hence
 
$$
P(A\cap \Omega,\Bbb{R}^n) \leq P(A,\Bbb{R}^n).\eqnlbl{3.4}
$$
 
Define $F' = \left (F\setminus B(x_0,\varepsilon ) \right)\cup (F\cap \overline{\Omega})$,
clearly
 
$$\aligned
F' \setminus  \overline{\Omega} & = \left( F \setminus  B(x_0,\varepsilon )\right) \setminus 
\overline{\Omega} \\
& = (F\setminus \overline{\Omega}) \setminus  B(x_0,\varepsilon ) \\
& = E_t \setminus  \overline{\Omega} \setminus  B(x_0,\varepsilon ) \\
& = \Cal{L}_t \setminus  \overline{\Omega} \setminus  B(x_0,\varepsilon ) \\
& = \Cal{L}_t \setminus  \overline{\Omega}
\endaligned
$$
and $F'\supset L_t$. Thus $F'$ is admissible in \eqnref{2.21} and therefore
$$
P(E_t,\Bbb{R}^n) \leq P(F',\Bbb{R}^n).
$$

Now we will show that $P(F',\Bbb{R}^n) \leq P(F,\Bbb{R}^n)$ which,
with the previous inequality, will imply \eqnref{3.3}.  First
observe from $E_t\Delta F\subset\subset B(x_0,\varepsilon )$ and
$(E_t\setminus \overline{\Omega}) \cap B(x_0,\varepsilon ) = (\Cal{L}_t \setminus 
\overline{\Omega}) \cap B(x_0,\varepsilon ) = \emptyset$ that $F'\cap
B(x_0,\varepsilon ) = F\cap B(x,\varepsilon ) \cap \overline{\Omega}$
and $F'\Delta F\subset\subset B(x_0,\varepsilon )$.  Hence we obtain
by \eqnref{3.4}
$$\aligned
P(F,\Bbb{R}^n) -  P(F',\Bbb{R}^n) & = P\left( F,B(x_0,\varepsilon )\right) -  P\left(
F',B(x_0,\varepsilon )\right) \\
& = P\left( F\cap B(x_0,\varepsilon ),B(x_0,\varepsilon )\right) -  P\left( F\cap
B(x_0,\varepsilon ) \cap \Omega, B(x_0,\varepsilon )\right) \\
& = P\left( F\cap B(x_0,\varepsilon ),\Bbb{R}^n\right) -  P\left( F\cap
B(x_0,\varepsilon ) \cap \Omega,\Bbb{R}^n\right) \\
& \geq 0
\endaligned\eqnlbl{3.5}
$$
This establishes \eqnref{3.3} when $g(x_0)<t$.

The argument to establish \eqnref{3.3} when $g(x_0)>t$ requires a
slightly different treatment from the previous case.  Since $G(x_0) =
g(x_0)>t$, the continuity of $G$ in $\Omega^c$ implies that
$\overline{B(x_0,\varepsilon )}\setminus \Omega\subset \Cal{L}_t$,
provided $\varepsilon $ is sufficiently small. Thus, we have
$\overline{B(x_0,\varepsilon )}\setminus \Omega\subset E_{t}$.
Clearly, we may assume $\varepsilon $ chosen to be smaller than
$\varepsilon _{0}$ of \eqnref{3.1}. Observe that the assumption that
$\partial \Omega $ is locally Lipschitz implies
that $P(\Omega ,B(x_{0},\varepsilon
))=P(\rn\setminus\Omega,B(x_{0},\varepsilon ))$. Consequently, we 
can appeal to \eqnref{3.1} to conlclude that $\rn\setminus\Omega $ is sub
area-minimizing in $B(x_{0},\varepsilon )$. On the other hand,
$E_{t}$ is super area-minimizing. Since $E_{t}\cap
B(x_{0},\varepsilon )\setminus\Omega \supset B(x_{0},\varepsilon
)\setminus\Omega $, we may apply Theorem
\thmref{maxprin} to find that $\partial
E_{t}=\partial(\rn\setminus\Omega )$ in some open neighborhood $U$ of $x_{0}$.
This implies that $L_{t}\cap U=\emptyset$ since $L_{t}\subset E_{t}$
and $\partial E_{t}\cap U\cap\Omega =\emptyset$. Consequently,
$E_{t}$ must be area-minimizing in $U$.

Thus far, we have shown that if either $g(x_{0})>t$ or $g(x_{0})<t$,
then $\partial E_{t}$ is area minimizing in a neighborhood of
$x_{0}$, say $B(x_{0},\varepsilon )$. We will show this leads to
a contradiction.  Assume first that $g(x_0) < t$ so that $G<t$ on
$(\Bbb{R}^n \setminus \Omega) \cap B(x_0,\varepsilon )$ provided that
$\varepsilon $ has been chosen sufficiently small.  Consequently
 
$$
E_t \cap B(x_0,\varepsilon ) \subset \overline{\Omega} \cap B(x_0,\varepsilon ).
\eqnlbl{3.7}
$$

We recall the notation concerning the representation of
$\partial\Omega$ as the graph of a Lipschitz function that preceeded
Lemma 3.1.  Thus with $x_0 \in \partial E_t \cap \partial \Omega\setminus
g^{-1}(t)$, we express $\partial\Omega$ locally around $x_0$ as
$\left\{ \left( x'h(x')\right):x'\in \right.$ $\left.
B(x'_0,\varepsilon )\right\}$, where $x_0 = (x'_0, x''_0)$ and $x''_0
= h(x''_0) > 0$.  For simplicity of notation, we take $x'_0 = 0$.
The number $\varepsilon '$ is chosen so that $\varepsilon
'<\varepsilon $ and

$$
\left\{ \left( x',h(x')\right):|x'| < \varepsilon '\right\} \subset B(x_0,\varepsilon ).
\eqnlbl{3.8}
$$
 
we define the half-infinite cylinder above $B'(0,\varepsilon ')$ as $C =
B'(0,\varepsilon ') \times [0,\infty]$.  Because of the local nature of the argument
we may assume that $$\Omega \cap C = \left\{ (x',x''):|x'|, 0\leq x'' <
h(x')\right\}.$$

Now consider the solution to the minimal surface equation on $B'(0,\varepsilon ')$
relative to the boundary data $f = h|_{\partial B'(0,\varepsilon ')}$ [M.M.
Chapter. 3].  Thus we let $v$ be the unique solution of
 
$$\aligned
div \left( \frac{\nabla v}{\sqrt{1+|\nabla v|^2}}\right) = 0 &\qquad
\text{on}\ B'(0,\varepsilon '). \\
v = f & \qquad \text{on}\ \partial B'(0,\varepsilon ')
\endaligned
$$
 
Since $h$ is a weak supersolution of the minimal surface equation, by Lemma
3.1, we have that $h\geq v$ on $\overline{B'(0,\varepsilon )}$, cf [GT, Theorem
10.7]. ; In fact, $h>v$ on $B'(0,\varepsilon ')$ because the set $\{h=v\}$ is
obviously closed in $B'(0,\varepsilon ')$ and it is also open in $B'(0,\varepsilon ')$
because of Lemma 3.2.  Hence, if this set is non-empty, $h=v$ in
$B'(0,\varepsilon ')$ which would contradict \eqnref{3.2}.  Consequently, with
$\delta = h(0) - v(0)$, we have $\delta > 0$.  Now consider a 1-parameter
family of graphs $v_\tau(x') = v(x') + \tau$ and let
 
$$
\tau_{m} = \max\left\{ \tau: \text{there exists $x''\in
\overline{B'(0,\varepsilon ')}$ such that : $\left(x',v_\tau(x')\right) \in
\partial E_t \cap \overline{\Omega}$}\right\}
$$
 
Note that $\tau_{m}\geq \delta$ since $v(0) + \delta = h(0)$ and $\left(
0,h(0)\right)\in \partial E_t \cap \overline{\Omega}$.  Let $V_{\tau_{m}} =
\left\{(x',x'') : |x'|<\varepsilon ', x''\leq v(x') + \tau_{m}\right\}$ and in view
of our choice of $\varepsilon '$, observe that
 
$$
E_t \cap \{x:|x'|<\varepsilon '\}\subset V_{\tau_{m}}.
$$
 
Observe also that if a point $\left( x',v_{\tau_{m}}(x')\right)$ is an element of
$\partial E_t \cap \overline{\Omega}$, then $|x'|<\varepsilon '$, for otherwise we
would have $v(x') + \tau_{m}\leq h(x')$ which would imply that $\tau_{m}\leq 0$. 
Thus the set $\partial\left[ E_t \cap \{x:|x'|<\varepsilon '\}\right] \cap \left\{
\left( x',v_{\tau_{m}}(x')\right) : |x'|<\varepsilon '\right\}$ is non-empty and
according to Theorem 2.2, it is open as well as closed in the connected set
$\left\{\left( x',v_{\tau_{m}}(x')\right) : |x'|<\varepsilon '\right\}$.  This
implies that
 
$$
\partial \left[ E_t\cap \{x:|x'|<\varepsilon '\}\right] \supset \left\{ \left(
x',v_{\tau_{m}}(x')\right) : |x'|<\varepsilon '\right\}. \eqnlbl{3.9}
$$
 
Since $t_{m}>0$, it follows that $v_{\tau_{m}}(x') > h(x'')$ whenever
$|x'|=\varepsilon '$.  consequently, using the continuity of $v_{\tau_{m}}$, the
graph $\left\{\left( x',v_{\tau_{m}}(x')\right),|x'|<\varepsilon '\right\}$ contains
points in $\Bbb{R}^n \setminus  \overline{\Omega}$, say $\left(
y',v_{\tau_{m}}(y')\right),|y'|<\varepsilon '$, as well as points in
$\overline{\Omega} \cap B(x_0,\varepsilon )$, say $\left( z',v_{\tau_{m}}(z')\right),
|z'|<\varepsilon '$.  The point $\left( y',v_{\tau_{m}}(y')\right),|y'|<\varepsilon '$,
could possibly be an element of $\Bbb{R}^n\setminus B(x_0,\varepsilon )$.  Consider the line
segment, $\lambda$, in $B'(x'_0,\varepsilon ')$ that joins $y'$ and $z'$.  Let $a'$
be the closet point to $y'$ on $\lambda$ with the property that $\left(
a',v_{\tau_{m}}(a')\right) \in \partial \Omega$.  Then all points $a$ on
$\lambda$ that are closer to $y'$ than $a'$ and that are sufficiently close to
$a'$ have the property that $\left( a,v_{\tau_{m}}(a)\right) \in (\Bbb{R}^n \setminus 
\overline{\Omega}) \cap B(x_0,\varepsilon )$.  Here we have used \eqnref{3.8} and
the continuity of $v_{\tau_{m}}$.  In view of \eqnref{3.9}, this implies that
$E_t\cap B(x_0,\varepsilon ) \cap (\Bbb{R}^n \setminus  \partial\Omega) \ne \emptyset$,
contradicting \eqnref{3.7}.  This contradiction was reached under the
assumption that $g(x_0) < t$ and the fact that $E_t$ is
area-minimizing in $B(x_0,\varepsilon )$.  A similar argument is employed in case
$g(x_0)>t$.
\hfill$\square$
\enddemo

In order to ultimately identify $E_t \cap \overline{\Omega}$ as the set $\{u
\geq t\}$ (up to a set of measure zero) for almost all $t$, we will need the
following result.

\proclaim{\thmlbl{goodcon} Lemma}  If $s,t\in T$ with $s<t$, then $E_t\subset\subset E_s$.
\endproclaim

\demo{Proof}  We first show that $E_t\subset E_s$.  Note that
 
$$\aligned
(E_s \cap E_t) \setminus  \overline{\Omega} & = (E_s \setminus  \overline{\Omega}) \cap (E_t \setminus 
\overline{\Omega}) \\
& = (\Cal{L}_s \setminus  \overline{\Omega}) \cap (\Cal{L}_t \setminus  \overline{\Omega}) \\
& = \Cal{L}_t \setminus  \overline{\Omega}
\endaligned
$$
and
$$
L_t \subset E_t, L_t \subset E_S \Longrightarrow L_t \subset E_s \cap E_t.
$$

Thus, $E_s \cap E_t$ is a competitor with $E_t$.

Similarly, 
 
$$\aligned
(E_s \cup E_t) \setminus  \overline{\Omega} & = (E_s \setminus  \overline{\Omega}) \cup (E_t \setminus 
\overline{\Omega}) \\
& = (\Cal{L}_s \setminus  \overline{\Omega}) \cup (\Cal{L}_t \setminus  \overline{\Omega}) \\
& = \Cal{L}_s \setminus  \overline{\Omega}
\endaligned
$$
 
and
 
$$
L_t \subset E_t, L_s \subset E_s \Rightarrow L_s \subset E_s \cup E_t.
$$
So $E_s \cup E_t$ is a competitor with $E_s$.
Then employing \eqnref{2.8}, we have
$$
P(E_{s})+P(E_{t})\leq P(E_s \cup E_t, \Bbb{R}^n) + P(E_s \cap E_t,\Bbb{R}^n) \leq P(E_s,\Bbb{R}^n) +
P(E_t,\Bbb{R}^n),
$$
and thus, since $E_{t}$ and $E_{s}$ are minimizers,
$$
P(E_s\cup E_t,\Bbb{R}^n) = P(E_s,\Bbb{R}^n)
$$
and
$$
P(E_s\cap E_t,\Bbb{R}^n) = P(E_t,\Bbb{R}^n).
$$
Reference to \eqnref{2.22} yields $|E_s\cup E_t| = |E_s|$, which in turn
implies $|E_s\setminus E_t| = 0$.  In view of \eqnref{2.9}
$$
x \in E\quad \text{if and only if
\quad$\limsup_{r\to 0} \frac{|E\cap B(x,r)|}{|B(x,r)|} > 0$}
$$
we conclude that $E_t \subset E_s$.

Now we come to the crucial part of the argument which is to show that
this containment is in fact strict.  For this 
purpose, first note that
$$
E_t \setminus  \overline{\Omega} = \Cal{L}_t \setminus  \overline{\Omega}\subset\subset \Cal{L}_s
\setminus  \overline{\Omega} = E_s \setminus  \overline{\Omega}. \eqnlbl{3.10}
$$
Now observe that \thmref{Lemma 3.3}implies
$$
\partial E_t \cap \partial E_s \cap \partial \Omega = \emptyset. \eqnlbl{3.11}
$$
In review of \eqnref{3.10} and \eqnref{3.11}, it remains to show that
$$
\partial E_t \cap \partial E_s \cap \Omega = \emptyset \eqnlbl{3.12}
$$
in order to establish the lemma.  For this purpose, let $S \equiv
\partial E_s \cap \partial E_t \cap \Omega$.  Then for $x_0 \in S$,
there are three possible cases with case (ii) being the central issue
of this paper. 
\medskip
\roster
\romitem For any $\varepsilon >0$,  ${L_{s}} \cap {L_{t}}\cap
\Omega\cap B(x_0,\varepsilon )$ is non-empty.
\romitem $x_{0}\in\overline{L_{s}}$ and 
$B(x_{0},\varepsilon )\cap L_{t}=\emptyset$ for some $\varepsilon >0$.
\romitem $B(x_0,\varepsilon )\cap L_t =
\emptyset=B(x_0,\varepsilon ) \cap L_s$ for some $\varepsilon>0$,
thus implying that both $\partial E_s$ and  
$\partial E_t$ are area-minimizing in $B(x_0,\varepsilon)$.
\endroster
Next, we will prove that above 3 cases are impossible, i.e. $S = \emptyset$,
which implies that $E_t \subset\subset E_s$.

For case (i), we can choose some sequence $\{y_n\} \subset L_s \cap L_t$, such
that $\lim_{n\to\infty} y_n = x_0$.  Since $\psi $ is continuous, we have
$\lim_{n\to\infty}\psi (y_n) = \psi (x_0) \geq t$.  Since $t>s$, there exists an
$\varepsilon >0$, such that $B(x_{0},\varepsilon )\subset E_{s}$
which conradicts the fact that $x\in\partial E_{s}$.

For case (ii), first observe that $E_{s}$ is super area-minimizing
and that $E_{t}$ is area-minimizing near $x_{0}$. Since $E_{t}\subset
E_{s}$, if follows from the maximum principle that $\partial E_{s}$
and $\partial E_{t}$ agree in a neighborhood of $x_{0}$.

For case (iii), since $\partial E_s$ and $\partial E_t$ are area
minimizing in $B(x_0,\varepsilon )$ and $E_t \subset E_s$, we apply
the maximum principle again to conclude that $\partial E_t$ and
$\partial E_s$ agree in a neighborhood of $x_0$.

Now combining above (i), (ii) and (iii), we conclude that $S$ is area-minimizing
and consists only of components of $\partial E_s$ that do not intersect
$\partial \Omega$.

Let $C$ be a component of $S$ and sing$\,\partial E_s$ the set of singular points
of $\partial E_s$. Intuitively, $C$  is a bounding cycle that is
area-minimizing, which is impossible. 
Our next step is to rigorously show that $C$ is empty. We divide the
proof of this into the following three parts.
		\medskip

{\bf Part 1.}  There exists an open set $V\subset \Bbb{R}^n$ such that
$\overline{V}\subset\Omega$ and $\partial V\supset S$.

For this purpose, we first find an open set $U$ such that
		\medskip

\itemitem{(i)}  $\partial U$ is an $(n-1)$-manifold with finitely many
components,
\itemitem{(ii)}  $\partial U \cap \partial E_s = \emptyset$,
\itemitem{(iii)}  $U\cap E_s = C$,
\itemitem{(iv)}  $\overline{U}\subset \Omega$,
\itemitem{(v)}  $U$ is connected.
		\medskip

\noindent To find such a set, consider a smooth approximation, $\rho $,
to the distance function $d(x) = d(x,C)$, i.e., let $\rho \in
C^\infty(\Bbb{R}^n \setminus  C)$ be such that
 
$$
k^{-1} d(x) \leq \rho (x) \leq k\, d(x)
$$
for all $x\in\Bbb{R}^n$, where $k$ is a positive number, cf [Z, Lemma 3.6.1]. 
Since $C$ is relatively open in $\overline{\Omega \cap \partial E_s}$, it
follows that $\partial \left\{ x:\rho (x) < r \right\}  \cap \partial
E_s = \emptyset$ for all small values of $r$.  Moreover, by Sard's Theorem and
the Implicit Function Theorem, $\rho ^{-1}(r)$ is a smooth
$(n-1)$-manifold for almost all values of $r$. For any such value of
$r$, let $U$ be the 
component of $\left\{ x:\rho (x)<r\right\}$ that contains $C$ to produce
a set satisfying all conditions (i)--(v) above except possibly (iv).  By
choosing $r$ sufficiently small, this too will be satisfied because $(\partial
E_s \cap \partial E_t) \cap \partial \Omega = \emptyset$.

Using only the fact that $\partial U$ is a compact $(n-1)$-manifold, we invoke
Alexander Duality of algebraic topology to conclude that $\Bbb{R}^n-\partial U$
consists of finitely many components, one more than the number of components in
$\partial U$, [GH, Theorem 27.10].  Moreover, each component of $\partial U$ is
the boundary of precisely one bounded open set.  Note that $\partial U_\infty$
is connected, where $U_\infty$ denotes the unbounded component of
$\Bbb{R}^n-\partial U$.  Indeed, since $U_\infty$ is connected, it is one of the
components of $\Bbb{R}^n-\partial U$.  thus there is a one-to-one
correspondence between the bounded components of $\Bbb{R}^n - \partial U$ and
the components of $\partial U$ which implies that $\partial U_\infty$ is
connected.

Since $\partial U_\infty$ is connected, either $\partial U_\infty \subset
(E_s)^i$ or $\partial U_\infty \subset (E_s)^e$, because $\partial U_\infty
\cap \partial E_s = \emptyset$.  In case $\partial U_\infty \cap (E_s)^i =
\emptyset$, define  $V$ by $V = (U_\infty)^e \cap (E_s)^i$.  Since $(U_\infty)^e
\supset U$ and $U\cap \partial E_s = C$, it follows that
$$
\partial V = \left( \partial U_\infty \cap (E_s)^i \right) \cup \left( (
U_\infty)^e \cap \partial E_s\right) \supset C.
$$
Similarly, if $(E_s)^e\cap \partial U_\infty = \emptyset$, define $V$ by $V =
(U_\infty)^e \cap (E_s)^e$, so that
$$
\partial V = \left( \partial U_\infty \cap (E_s)^e\right) \cup \left( (
U_\infty)^e \cap \partial E_s\right)\supset C.
$$
Thus we have established the existence of an open set $V$ that is either a
subset of $(E_s)^i$ or a subset of $(E_s)^e$ and satisfies
$$
\overline{V} \subset \Omega, \partial V \supset C.
$$

To finish the proof of the Lemma, we will now show that this leads to a
contradiction.
		\medskip

{\bf Part 2.}  If $V \subset (E_s)^i$, then $E_s$ is not a minimizer of
\eqnref{2.21}.

There are two cases here.

Case 1.  $V\cap L_s = \emptyset$.  This implies that the closed set $F_s$
defined by $F_s = \overline{E_s - V}$ is admissible in the minimization problem
\eqnref{2.21}.

If we can show that
 
$$
H^{n-1} (\Omega \cap \partial E_s) \geq H^{n-1} (\Omega \cap \partial F_s) +
H^{n-1} (\partial V), \eqnlbl{3.17}
$$
 
the desired conclusion is reached since then $H^{n-1}(\Omega \cap \partial E_s)
> H^{n-1} (\Omega \cap \partial F_s)$, contradicting the minimality of
$H^{n-1}(\Omega\cap \partial E_s)$.

To establish \eqnref{3.17}, it is sufficient to prove
 
$$
\partial F_s \cap \partial V = \emptyset\eqnlbl{3.18}
$$
since $\Omega \cap \partial E_s \supset \left[ ( \Omega \cap \partial
F_s) - (\partial F_s) \cap (\partial V)\right] \cup \partial V$.
Because $\partial V = (U_\infty)^e \cap \partial E_s$, it follows
that for all sufficiently small $r>0$,
$$
B(x,r) \cap \partial V = B(x,r) \cap \partial E_s.
$$
Furthermore, for all small $r>0$
$$
B(x,r) \cap \overline{V} = B(x,r) \cap E_s.
$$
 
It follows immediately that $x\not\in \partial F_s$ and therefore \eqnref{3.18}
is established.
\medskip

Case 2.  $V\cap L_s \ne \emptyset$, but $V\cap L_t = \emptyset$.  Then replace
$E_s$ and $F_s$ by $F_t = \overline{E_t \setminus  \overline{V}}$ as in Case 1, follow
the same line of proof as in Case 1, we will conclude that $E_t$ is not a
minimizer of \eqnref{2.21}.
		\medskip

Case 3.  $V\cap L_s \ne \emptyset,\ V\cap L_t \ne \emptyset$.  Then there exist
components $L_{s_1}$ and $L_{t_1}$ of $L_s$ and $L_t$ respectively such that
$L_{t_1}\subset\subset L_{s_1}\subset V$. Since $\partial
V\subset\partial E_{s}$, it follows that $V$ is a minimizer of the of
the obstacle problem with $L_{t_{1}}$ as the obstacle. But this is
not possible since $L_{t_1}\subset\subset V$.
		\medskip

{\bf Part 3.}  If $V\subset (E_s)^e$, then $E_s$ is not a minimizer of
\eqnref{2.21}.  Let $G_s = E_s \cup \overline{V}$, then $G_s$ is an admissible
competitor in \eqnref{2.21}.  Now repeat the argument of Part 2 case 1 with
$F_s$ replaced by $G_s$ to contradict the minimality of $\partial E_s$.
\qqed
\enddemo

We now are in a position to construct the solution $u$ to problem
\eqnref{1.3}. For this purpose, we first define for $t\in T$,
$$
A_{t}=\overline{E_{t}\cap\Omega}.
$$
With the help of Lemma \thmref{goodcon}, observe that for $t\in
T$, 
$$
\align
\{g>t\}&\subset (E_{t})^{i}\cap\partial \Omega\subset
A_{t}\cap\partial \Omega\eqnlbl{gw}\\
\overline{\{g>t\}}&\subset
A_{t}\cap\partial \Omega\subset\overline{E_{t}}\cap\partial \Omega
=[(E_{t})^{i}\cup\partial E_{t}]\cap\partial \Omega\subset\{g\geq
t\}.\eqnlbl{lvlset}
\endalign 				
$$
Finally, note that \eqnref{lvlset} and Lemma \thmref{scontain} imply
$$
A_{t}\subset\subset A_{s}\eqnlbl{acontain}
$$
relative to the topology on $\obar$ whenever $s,t\in T$ with $s<t$.
We now define our solution $u$ by 
$$
u(x)=\sup\{t:x\in A_{t}\}.\eqnlbl{defu}
$$
\bigskip\noindent
\proclaim{\thmlbl{dsoln} Theorem} The function $u$ defined by
\eqnref{defu} satisfies the following:
\itemitem{\rm(i)} $u=g$ on $\partial \Omega$
\itemitem{\rm(ii)} $u$ is continuous on $\obar$,
\itemitem{\rm(iii)} $A_{t}\subset\{u\geq t\}$ for all $t\in T$ and
$\abs{\{u\geq t\}-A_{t}}=0$ for
almost all $t\in T$.
\itemitem{\rm(iv)} $u\geq\psi $ on $\overline\Omega $.
\endproclaim
\demo{Proof} To show that $u=g$ on $\partial \Omega$, 
let $x_{0}\in\partial \Omega$ and suppose $g(x_{0})=t$. If $s<t$, 
then $G(x)>s$
for all $x\in\Omega^{c}$ near $x_{0}$. Hence, $x_{0}\in
(E_{s})^{i}\cap\partial 
\Omega$ by \eqnref{gw} and consequently, $x_{0}\in A_{s}$
for all $s\in T$ such that $s<t$. By \eqnref{defu}, 
this implies $u(x)\geq t$. To show that $u(x)=t$ suppose by
contradiction that $u(x)=\tau >t$. Select $r\in(t,\tau)\cap T$. Then $x\in
A_{r}$. But $A_{r}\cap\partial \Omega\subset\{g\geq r\}$ by
\eqnref{lvlset}, a 
contradiction since $g(x)=t<r$.

For the proof of (ii), it is easy to verify that
$$
\{u\geq t\}=\left\{\medcap A_{s}:s\in
T,s<t\right\}\;\hbox{and}\;\{u>t\}=\left\{\medcup A_{s}:s\in T,s>t\right\}. 
$$
The first set is obviously closed while the second is open relative
to $\obar$ by
\eqnref{acontain}. Hence, $u$ is continuous on $\obar$.

For (iii), it is clear that $\{u\geq t\}\supset A_{t}$.
Now, $\{u\geq t\}-A_{t}\subset u^{-1}(t)$. But
$\abs{u^{-1}(t)}=0$ for almost all $t$ because
$\abs{\Omega}<\infty$.

In (iv), it is sufficient to show $u(x_{0})\geq\psi (x_{0})$ for
$x_{0}\in\Omega $. Let $t=u(x_{0})$ and $r=\psi (x_{0})$ and suppose
$t<r$. Then $x_{0}\in L_{r'}\subset E_{r'}$ for $t<r'<r$. But then,
$x_{0}\notin A_{r'}$ by the definition of $u$, a contradiction.
\qqed
\enddemo

\proclaim{\thmlbl{soln} Theorem} If $\Omega$ is a bounded Lipschitz
domain that satisfies \eqnref{3.1}
and \eqnref{3.2}, then the function $u$ defined by
\eqnref{defu} is a solution to \eqnref{1.3}.
\endproclaim
\demo{Proof} Let $v\in BV(\Omega),\;v=g$ on $\partial \Omega$ be a
competitor in problem \eqnref{1.3}. We recall the extension $G\in
BV(R^{n}-\obar)$ of $g$, \eqnref{2.17}. Now define an extension
$\vbar\in BV(R^{n})$ of $v$ by $\vbar=G$ in $R^{n}-\obar.$ Let
$F_{t}=\{\vbar\geq t\}$. It is sufficient to show that 
$$
P(E_{t},\Omega)\leq P(F_{t},\Omega)\eqnlbl{min}
$$
for almost every $t\in T$ (see \eqnref{2.20}), because then $v\in BV(\Omega)$ and
\eqnref{coarea} would imply 
$$
\int_{a}^{b}P(E_{t},\Omega)\;dt
\leq\int_{-\infty}^{\infty}P(F_{t},\Omega)\;dt=\norm{\nabla v}(\Omega)<\infty.
$$
Hence, by \eqnref{reversecoa}, $u\in BV(\Omega)$; furthermore,
$\norm{\nabla u}(\Omega)\leq\norm{\nabla v}(\Omega)$ by \eqnref{coarea}.

We know that $E_{t}$ is a solution of 
$$
\min\left\{ P(E,\Bbb{R}^n): E \setminus  \overline{\Omega} = \Cal{L}_t \setminus 
\overline{\Omega}, \;E\supset L_t\right\}, \eqnlbl{2.22}
$$
while $F_{t}-\obar=\l{t}-\obar$ and $F_{t}\supset L_{t}$. Hence, 
$$
P(E_{t},R^{n})\leq P(F_{t},R^{n}).\eqnlbl{q}
$$
Next, note that
$$
\eqalign{
P(E_{t},R^{n})&=H^{n-1}(\partial _{*}E_{t}-\obar)+H^{n-1}(\partial
_{*}E_{t}\cap\partial \Omega)+H^{n-1}(\partial _{*}E_{t}\cap\Omega)\cr
&\geq H^{n-1}(\partial _{*}\l{t}-\obar)+P(E_{t},\Omega).\cr}\eqnlbl{change}
$$
We will now show that
$$
\eqalign{
P(F_{t},R^{n})&=H^{n-1}(\partial _{*}\l{t}-\obar)+H^{n-1}(\partial
_{*}F_{t}\cap\Omega).\cr 
&=H^{n-1}(\partial _{*}\l{t}-\obar)+P(F_{t},\Omega),\cr}\eqnlbl{star}
$$
which will establish \eqnref{min} in light of \eqnref{q} and \eqnref{change}.

Observe 
$$
P(F_{t},R^{n})=H^{n-1}(\partial _{*}\l{t}-\obar)+H^{n-1}(\partial
_{*}F_{t}\cap\partial \Omega)+H^{n-1}(\partial _{*}F_{t}\cap\Omega).
$$
We claim that $H^{n-1}(\partial _{*}F_{t}\cap\partial \Omega)=0$ for
almost all $t$ because $\partial _{*}F_{t}\subset\partial
F_{t}\subset\overline{v}^{-1}(t)$ since $\overline{v}\in
C^{0}(R^{n})$. But $H^{n-1}(\overline{v}^{-1}(t)\cap\partial
\Omega)=0$ for all but countably many $t$ since $H^{n-1}(\partial
\Omega)<\infty$. 
\qqed
\enddemo

\chapno=5
\tagno=0
\thmno=0
\subsectionnumber=0
\Increase{\sectionnumber}
\head
\the\chapno. Modulus of continuity of the solution
\endhead
\proclaim
{\thmlbl{holderatbdry} Lemma} Suppose $\Omega$ is a bounded,
open subset of $\R^{n}$ whose boundary is $C^{2}$ with mean curvature
bounded below by $a>0$. Assume $g\in C^{0,\alpha}(\partial \Omega)$,
and $\psi \in C^{0,\alpha /2}(\Omega )$ 
for some $0<\alpha\leq1$. Let $u\in C^{0}(\obar)\cap BV(\Omega)$ be a solution
to \eqnref{3.1}. Then, there exist positive numbers $\delta$ and
$C$ depending only on $a, \norm{g}_{C^{0,\alpha}(\partial \Omega)},\hfill\break
\norm{g}_{C^{0}(\partial \Omega )}, \norm{\psi }_{C^{0,\alpha
/2}(\Omega )}$ and $\norm{u}_{C^{0}(\Omega) }$ such that
$$
\abs{u(x)-u(x_{0})}\leq C\abs{x-x_{0}}^{\alpha /2}.
$$
wherever $x_{0}\in\partial \Omega $  and $x\in\overline{\Omega }$ with
$\abs{x-x_{0}}<\delta$.
\endproclaim
\demo{Proof} For each $x_{0}\in\partial \Omega$ we will constuct
functions $\wp,\wm\in C^{0}(\overline{U})$ where\hfill\break
$U(x_{0},\delta ):=B(x_{0},\delta)\cap\Omega$ and $\delta >0$ is sufficiently small, such that 

\roster
\romitem $\wp(x_{0})=\wm(x_{0})=g(x_{0})$,
\romitem For $x\in U(x_{0},\delta )$ 
$$
\align
\abs{\wp(x)-g(x_{0})}&\leq
C\abs{x-x_{0}}^{\alpha/2}\\
\abs{\wm(x)-g(x_{0})}&\leq
C\abs{x-x_{0}}^{\alpha/2},
\endalign
$$
\romitem $\wm\leq u\leq\wp$ in $U(x_{0},\delta )$
\endroster

We begin with the construction of $\wm$. To
this end, let $d(x)=\,\text{dist}\,(x,\partial \Omega)$. Since
$\partial \Omega\in C^{2}$ recall that $d\in C^{2}(\{x:0\leq
d(x)<\delta_{0}\})$ for some $\delta_{0}>0$. Furthermore, since 
$\partial \Omega$ has positive mean
curvature and $\abs{\nabla d}=1$, it follows that 
$$
\div\left(\frac{\nabla d}{\abs{\nabla d}}\right)=\Delta d\leq -a, \eqnlbl{lap}
$$ 
in $\{x:0\leq d(x)<\delta_{0}\}$ for some $a>0$. For each
$\varepsilon >0$, set
$$
\eqalign{
v(x)&=\abs{x-x_{0}}^{2}+\lambda d(x)\cr
\wm(x)&=\max\{\psi, -Kv^{\alpha/2}(x)+g(x_{0})\},\cr}
$$
where $\lambda>0$ is to be determined later. Clearly (i) is satisfied.

Next, in the open set
$\{\wm>\psi \}$, observe that
$$
\eqalign{
\abs{\nabla
\wm}&=K{\textstyle{\frac{\alpha}{2}}}v^{\frac{\alpha}{2}-1}\abs{\nabla 
v},\cr 
\abs{\nabla v}&=\abs{2(x-x_{0})+\lambda\nabla d}\geq\lambda\abs{\nabla
d}-2\abs{x-x_{0}}\cr
&=\lambda-2\abs{x-x_{0}}>0,\cr}
$$
provided we choose $\delta$ and $\lambda$ such that
$\lambda>2\delta$. Further, we note that 
$$
\div\left(\frac{\nabla \wm}{\abs{\nabla \wm}}\right)=
-\div\left(\frac{\nabla v}{\abs{\nabla v}}\right)=\frac{-1}{\abs{\nabla
v}^{3}}Av,
$$
where $Av=\abs{\nabla v}^{2}\Delta v-D_{i}vD_{j}vD_{ij}v$. Finally,
observe that $Av<0$ for $\lambda$ sufficiently large and $\delta$
sufficiently small. Indeed, using $D_{i}dD_{ij}d=0$ for any $j$, one
readily obtains
$$
Av=\abs{\nabla v}^{2}(\lambda\Delta d+2(n-1))-4\lambda(x-x_{0})_{i}(x-x_{0})_{j}D_{ij}d
$$
and
$$
\abs{\nabla v}^{2}=\lambda^{2}+4\abs{x-x_{0}}^{2}+4\lambda(x-x_{0})\cdot\nabla d
$$
so that 
$$
Av\leq-a \lambda^{3}+O(\lambda^{2}),\;\text{\rm as }\;\lambda\to\infty
$$
uniformly for $\delta<\delta_{0}$.  

Clearly, we can choose $K$ sufficiently large so that
$\wm=\psi$ on
$\partial U(x_{0},\delta )$ and that (ii) is satisfied, where $K$ depends only on
$\norm{g}_{C^{0,\alpha}(\partial \Omega)}, \norm{g}_{C^{0}(\partial \Omega
)}, \norm{\psi }_{C^{0,\alpha /2}(\Omega )}$ and $\norm{u}_{C^{0}(\Omega) }$.
Also, on $\Delta :=\{\wm>u \}\cap U(x_{0},\delta )$, we have 
$\wm=-Kv^{\alpha /2}+g(x_{0})$ and 
therefore
$$
\abs{\nabla \wm}>0\quad\text{and}\quad\div\left(\frac{\nabla
\wm}{\abs{\nabla \wm}}\right)>0\;\text{on $\Delta $}.\eqnlbl{posmeancurvature}
$$

We now proceed to show that $\Delta =\emptyset$, which will
establish the first of the inequalities in (iii). For this purpose,
note that $\wm\in BV(\Delta )$. Next, for $t>0$, let $\Delta
_{t}:=\{\wm-t>u\}$ and note that 
$$
\Delta =\medcup_{t>0}\Delta _{t}\qquad \Delta
_{t}\subset\subset\Delta\subset\Omega.\eqnlbl{propercontain} 
$$
Let $\omega ^{*}:=\max(u,\wm-t)$ and note that $\omega
^{*}\in BV(\Omega )\cap C^{0}(\overline\Omega )$ since $\wm-t=\psi
-t<u$ on $\partial \Delta _{t}$. 
For all but countably many $t>0$, it follows from basic measure theory that
$$
\norm{\nabla \omega ^{*}}(\partial \Delta _{t})=0=\norm{\nabla
u}(\partial \Delta _{t}).\eqnlbl{goodt}
$$
For the remainder of this argument, we will consider only such $t$.
Since $\omega ^{*}\geq u\geq\psi $, it follows
that 
$$
\norm{\nabla u}(\Omega )\leq\norm{\nabla \omega ^{*}}(\Omega ).\eqnlbl{badinequality}
$$

Now let $\eta\in C^{\infty}_{0}(\Delta )$ satisfy $\eta=1$ on $\Delta _{t}$ and
$0\leq\eta\leq1$ in $\Delta $. Set
$$
h=\eta\frac{\nabla \wm }{\abs{\nabla \wm }}
$$
so that $h\in[C^{1}_{0}(\Delta )]^{n}$. Since $\omega ^{*}=u$ on
$\Delta -\Delta_{t}$, it follows from
$$
\leqalignno{
\int_{\Delta }u\;\div h\;dx&=-\nabla u(h),\cr
\noalign{and}
\int_{\Delta }\omega ^{*}\;\div h\;dx&=-\nabla \omega ^{*}(h),&\cr}
$$
that 
$$
\int_{\Delta }u-\omega ^{*}\;dx=\int_{\Delta_{t}}(u-\wm +t )\;\div h\;dx=[\nabla (\omega ^{*}-u)](h).
$$
It follows from \eqnref{goodt}and the definition of the BV norm 
that 
$$
\norm{\nabla \omega ^{*}}(\partial \Delta_{t})\leq\norm{\nabla u}(\partial
\Delta_{t})+\int_{\partial \Delta_{t}}\abs{\nabla \omega ^{*} }\;dx=0,
$$
so that
$$
\eqalign{
\int_{\Delta_{t}}(u-\omega ^{*} +t )\;\div h\;dx&=\nabla \omega ^{*} (h\upchi_{\Delta_{t}})-\nabla
u(h\upchi_{\Delta_{t}})\cr
&\geq\int_{\Delta_{t}}\abs{\nabla \omega ^{*} }\;dx-\norm{\nabla u}(\Delta_{t}).\cr} 
$$
Since $u-\wm +t <0$ and $\;\div h>0$ on $\Delta_{t}$, we have
$$
\int_{\Delta_{t}}\abs{\nabla \omega ^{*} }\;dx<\norm{\nabla u}(\Delta_{t}).
$$
That is, 
$$
\norm{\nabla \omega ^{*}}(\Delta_{t} )<\norm{\nabla u}(\Delta_{t} ).
$$
Since $\omega ^{*}=u$ on $\rn\setminus\Delta _{t}$, we obtain 
from \eqnref{goodt} that $\norm{\nabla \omega ^{*}}(\Omega )<\norm{\nabla u}(\Omega )$,
which contradicts \eqnref{badinequality}. Thus we conclude that
$\wm\leq u$ on $U(x_{0},\delta )$. 

The proof of the second inequality in (iii) is obtained by a similar
argument using $\wp(x):=Kv^{\alpha /2}(x)+g(x_{0}$.
\qqed\enddemo
\proclaim
{\thmlbl{holder} Theorem} Suppose $\Omega$ is a bounded,
open subset of $\R^{n}$ with $C^{2}$ boundary having mean curvature
bounded below by $a>0$. Suppose $g\in C^{0,\alpha}(\partial \Omega)$,
and $\psi \in C^{0,\alpha /2}$ 
for some
$0<\alpha\leq1$. If $u\in C^{0}(\obar)\cap BV(\Omega)$ is a solution
to \eqnref{3.1}, then $u\in C^{0,\alpha/2}(\obar)$.
\endproclaim
\demo{Proof}
For $s<t$, consider the superlevel sets $E_{s},E_{t}$ of $u$ and assume
that\hfill\break $\text{dist}(\partial E_{s},\partial E_{t})=\abs{y-x}$
where $x\in E_{t}$ and $x\in E_{s}$. Assume $t-s$ is small enough to
ensure that $\abs{y-x}<\delta$, where $\delta $ is given by Lemma
\thmref{holderatbdry}. Observe that $L_{t}\subset
E_{t}\subset\subset E_{s}$. Theorems \thmref{dsoln} and \thmref{soln}
imply that $u$ is continuous on $\overline{\Omega }$ and therefore
bounded. Hence it is sufficient to show that $\abs{u(y)-u(x)}=\abs{t-x}\leq
C\abs{x-y}^{\alpha /2}$ whenever $\abs{y-x}<\delta $. This will be
accomplished by examining the following five cases.

(i) If either $x$ or $y$ belongs to $\partial \Omega $, then our
result follows from Lemma \thmref{holder}.
(ii) $y\in\partial E_{s}\setminus\overline{L_{s}}$ and
$x\in\partial E_{t}\setminus\overline{L_{t}} $: Let
$[\partial E_{t}]_{v}$ denote the translation of $\partial E_{t}$ by
the vector $v:=(y-x)\abs{y-x}$. Since both $\partial E_{s}$ and
$\partial E_{t}$ are area-minimizing in some neighborhoods of $y$ and
$x$ respectively, we can apply Theorem \thmref{maxprin} to conclude
that $\partial E_{s}$ and $[\partial E_{t}]_{v}$ agree on some
connected component of $\partial E_{s}$, say $S$, that contains $y$.

If $S\cap\partial \Omega =\emptyset$, it follows that $S$ is
area-minimizing in some open set $U\supset S$. 
Now let $S'$ be a component of the set of regular points of
$S$. We first show that $S'$ is a cycle in the sense of currents;
that if, we wish to show that
$$
\int_{S'}\text{\rm d}\varphi=0\eqnlbl{cycle}
$$
whenever $\varphi$ is a smooth $(n-2)$-form supported in $B(0,R)$
where $B(0,R)$ is the ball having the property
that $\Omega\subset\subset B(0,R)$. Since $S'$ is area-minimizing in
$\Omega$, we appeal to the monotonicity formula [S1, \S17.6] to
conclude that only a finite number of components of $(\partial
E_{s})$ can intersect any given compact subset of $\Omega$, in
particular, $\text{\rm spt}\;\varphi\cap\overline{S'}$. Thus, there exists a
smooth function $\zeta$ that is $1$ on
$\text{\rm spt}\;\varphi\cap\overline{S'}$ and $0$ in a neighborhood of
$\partial E_{s}-\overline{S'}$. Then, \eqnref{cycle} is established by 
$$
\int_{S'}\text{\rm d}\varphi=\int_{S'}\text{\rm d}(\zeta
\varphi)=\int_{\partial E_{s}}\text{\rm d}(\zeta
\varphi)=0.
$$
Thus, $S'$ is an
$(n-1)$-rectifiable cycle in the sense of currents; that is,
$\partial S'=0$. Now appeal to [S1, 27.6] to conclude that
there is a measurable set $F\subset B(0,R)$ such that $\partial
F=S'$. It follows from elementary considerations that for a given
vector $\nu\in R^{n}$, there is a hyperplane, $P$, with normal $\nu$
such that $P\cap \overline{S'}\not=\nullset$ and 
$$
F\subset\{x:(x-x_{0})\cdot\nu\leq0\}
$$
where $x_{0}\in P\cap \overline{S'}$. Theorem \thmref{maxprin} implies
$P\cap \overline{S'}$  
is open as well as closed in $P$, thus leading to a contradiction in
case $S\cap\partial \Omega =\emptyset$.
If $S\cap\partial \Omega \neq\emptyset$, we then are led to a
situation covered by (i).

(iii) $y\in\partial E_{s}\setminus\overline{L_{s}},x\in\partial
E_{t}\cap\overline{L_{t}}$: Since $L_{s}\supset L_{t}$, there exists
$y'\in\partial L_{s}$ such that $y'-x= c(y-x), 0<c<1$, and therefore
$$
\abs{u(y)-u(x)}=\abs{t-s}=\abs{\psi (y')-\psi (x)}\leq
C\abs{y'-x}^{\alpha }\leq C\abs{y-x}^{\alpha }.\eqnlbl{psiholder}
$$

(iv) $y\in\partial E_{s}\cap\overline{L_{s}},x\in\partial
E_{t}\cap\overline{L_{t}}$: This is treated as in \eqnref{psiholder}.

(v) $y\in\partial E_{s}\cap\overline{L_{s}},x\in\partial
E_{t}\setminus\overline{L_{t}}$: In this case we can apply Corollary
\thmref{maxprinsuper} to obtain an area minimizing connected component $S$ of
$\partial E_{s}$ which can be treated as in (ii) above.
\qqed
\enddemo

\chapno=6
\subsectionnumber=0
\Increase{\sectionnumber}
\head
\the\chapno. A Monotonicity Principle for Superminimizing Sets
\endhead
\thmno=0
\tagno=0 

An issue left open in our development is whether
the regularity requirement $\overline{E}\cap U= \overline{E^i}\cap U$
is necessary in Theorem \thmref{maxprin}, the extended maximum
principle for sub and superminimizing sets.  

This suggests the question, of interest in its own right, of
what regularity, if any, is enjoyed by (sub)superminimizing sets.
For example, do (sub)superminimizers have tangent cones?  Are
they $C^1$ or analytic $H^{n-1}$ almost-everywhere?  And,
the question begged by Theorem \thmref{maxprin}, 
is a subminimizer necessarily the closure of its interior?
In the next section, we will give an explicit example showing that the
last conjecture is false.
In this section, we present some preliminary results in the
direction of regularity, consisting of a new {\it monotonicity principle}
and consequent {\it one-sided mass bound} for (sub)superminimizing sets.

Let $B_r=B(0,r)$ denote the ball of radius $r$ about the origin
in $\Bbb{R}^n$.
Let $F$ be a superminimizing set in $U$, and without loss of
generality, assume $B_1\subset U$.

\proclaim{\thmlbl{1} Lemma} 
Let $\tilde A=\{x\in A^c:\text{\rm the metric density of $A$ is one at
$x$}\}$. Then,
$ H^{n-1}(\partial B_r \cap \tilde A)=0$ for almost
all $r$.
\endproclaim

\demo{Proof}
The Lebesgue measure of $\tilde A\cap B_1$ is zero.
But, by the co-area formula, \eqnref{coarea}, it is also equal to
$\int_0^1 H^{n-1}(\partial B_r \cap \tilde A) \,dr$, whence
the result follows.
\qqed
\enddemo

\proclaim{\thmlbl{2} Lemma}  
Let $E$ area subminimizing in $U$, $B_1\subset U$, and
$r$ such that $H^{n-1}(\partial B_r \cap \tilde E)=0$.
Then, $P(E,B_r) \leq H^{n-1}(E\cap \partial B_r)$.
\endproclaim

\demo{Proof}
The set $G=E\setminus B_r$ is a 
competitor to $E$.  Exterior to $B_r$,
$G$ has the same reduced boundary as does $E$,
but interior to $B_r$, it has no reduced boundary.
On $\partial B_r$, $G$ has reduced boundary contained
in the set of points at which $E$ has density one,  
which by assumption is contained in $E$
except for a set of $H^{n-1}$-measure zero.

Therefore, by the subminimality of $E$, we have
$$
0\leq P(E \setminus B_r,U) - P(E,U)
\leq H^{n-1}(E \cap \partial B_r) -
P(E,B_r),
$$
giving the result.
\qqed
\enddemo

Define the dimension-dependent constant $0<\delta(n)<1/2$ by
$$
\delta(n)=|D_1|/|B_1|,
$$
where $D_1\subset B_1$ is a set bounded by a hemispherical cap of radius
one,  orthogonal to $\partial B_1$.


\proclaim{\thmlbl{3} Lemma}
If $|A\cap B_r|/|B_r| \leq \delta(n)$, then
$P(A,B_r)/H^{n-1}(\partial B_r) \geq |A\cap B_r|/|B_r|$.
\endproclaim

\bigskip
\noindent
{\bf Remark.}  Another way of stating this result is
that $P(A,B_r) \geq (n/r)|A\cap B_r|$.  It could also be
rephrased as an isoperimetric inequality.

\demo{Proof}
By rearrangement, we find that the set $D$ of minimum
perimeter $P(D,B_r)$ subject to $|D\cap B_r|=|A\cap B_r|$ is 
the set bounded by a hemispherical cap meeting $\partial B$ orthogonally.
Trivially, we have
$$
P(D,B)\leq P(A,B).
\eqnlbl{2}
$$

Let $D_r$ be the set bounded by a spherical cap of radius $r$,
intersecting $\partial B_r$ orthogonally, so that $|D\cap B_r|/|B_r|=\delta(n)$.
Since $|D\cap B_r|/|B_r|=|A\cap B_r|/|B_r| \leq \delta(n)$, we thus have that 
$|D|\leq |D_r|$ and so the radius of the hemispherical cap bounding $D$ is 
less than or equal to $r$.  It follows by elementary geometry that
$$
H^{n-1}(\partial D \cap \partial B) \leq P(D,B).
\eqnlbl{3}
$$
(To see this, e.g.,  one can reflect the hemispherical cap $D$
about the plane of its intersection with $B_r$, to obtain a
surface oriented in the same direction as the patch $\overline{D} \cap \partial B_r$
and containing  the patch in its interior.  Since the patch has positive
mean curvature, it follows that this outer surface has greater 
area than does  $\overline{D} \cap \partial B_r$.)   

But, $D$ is entirely contained in the cone $C$ from 
$\partial D \cap \partial B_r$ to the center of $B_r$
and tangent to $D$ at $\partial B_r$.
That is,   $|A\cap B_r|\leq |C|$.  On the other hand,
the volume ratio $|C|/|B_r|$ for a cone is exactly
its surface ratio, $H^{n-1}(\partial D \cap \partial B_r)/H^{n-1}(\partial B_r) $.
Combining these facts with \eqnref{3} and \eqnref{2}, we have
$$
\eqalign{
|A|/|B_r| &\leq |C|/|B_r| =  H^{n-1}(\partial D \cap \partial B_r)
/H^{n-1}(\partial B_r) \cr
&\leq
P(D,B_r) /H^{n-1}(\partial B_r)
\leq 
P(A,B_r) /H^{n-1}(\partial B_r), \cr
}
$$
which leads to our desired conclusion.
\qqed
\enddemo

We now prove our main result, a {\it volume} monotonicity principle
for superminimizing sets.

\proclaim{\thmlbl{4} Proposition}
Let $E$ be subminimizing in $U$, $B_1\subset U$.
If $|E\cap B_1|/|B_1| < \delta(n)$ ( $0<\delta(n)<1/2$ as
defined above \thmref{3}),  then the ratio
$|E\cap B_r|/|B_r|$ is increasing in $r$
for $0\leq r\leq 1$.
\endproclaim

\demo{Proof}
From Lemmas \thmref{2} and \thmref{3}, we have 
$$
H^{n-1}(E\cap \partial B_r)) /H^{n-1}(\partial B_r) \geq
P(E,B_r)/H^{n-1}(\partial B_r)  \geq 
|E\cap B_r|/|B_r)|$$
for almost all $r$, so long as 
$|E\cap B_r|/|B_r| < \delta(n)$.

By the co-area formula, \eqnref{coarea}, 
$$
(d/dr)|B_r|= H^{n-1}(\partial B)\; \text{and}\;
(d/dr)|E|= H^{n-1}(E \cap \partial B_r).
$$
Thus, 
$$
d|E\cap B_r|/d|B_r|
=
H^{n-1}(E \cap \partial B_r)/ H^{n-1}(\partial B)
\geq
|E\cap B_r|/|B_r|,$$
giving monotonicity so long as
$|E\cap B_r|/|B_r|<\delta(n)$.
But, because of monotonicity, this
property persists for all $0\leq r \leq 1$.
\qqed
\enddemo

This property has many implications.  Among them
is the following important one, a one-sided bound on
the average density.

\proclaim{Proposition \thmlbl{5}}
Let $E$ be subminimizing in $U$, $B_1\subset U$.
If $0\in \partial E$, then $|E\cap B_1|/|B_1|\geq  \delta(n)$.
\endproclaim

\demo{Proof}
Suppose to the contrary that $|E\cap B_1|/|B_1|<\delta(n)$.
Then, for some $R<1$,
$|E\cap B(x,R)|<\delta(n)$ for every $x\in B_{1-R}$.
By the monotonicity property of Proposition \thmref{4},
we thus have $|E\cap B(x,r)|/|B(x,r)| < \delta(n)$
for $r\leq R$.  Thus,
$$
|E\cap \tilde B|/|\tilde B| < \delta(n)<1/2
$$
for any ball contained in $B_{1-R}$; hence the density
of $E$ is strictly less than $1/2$ at each point of $B_{1-R}$.

But, since the density of $E$ must be zero or one at almost every
point of $B_{1-R}$, the density of $E$ must be zero at almost 
every point in $B_{1-R}$, and therefore $|E\cap B_{1-R}|=0$.
But, by our convention in choosing set representatives,
this would imply that $B_{1-R}\subset E^e$, in particular $0\in E^e$,
a contradiction.
\qqed
\enddemo

\proclaim{\thmlbl{5.1} Corollary} 
If $E$ is subminimizing,
then $\overline{E^i_m}= (E^e_m)^c=\overline{E}$.
\endproclaim

\demo{Proof}
By Proposition \thmref{5}, the density of $E$ at any $x\in \partial E$
is strictly greater than $0$, hence $\partial E\cap E^e_m=\emptyset$.
It follows that $\partial E$, and therefore $\overline{E}$ as well,
is contained in 
$ (E^e_m)^c \subset \overline{E_m^i}$.
Since $\overline{E_m^i}$ is always contained in $\bar{E}$, we thus 
obtain
$$
\bar{E^i_m}= (E^e_m)^c =\bar{E},
$$
as claimed.
\qqed
\enddemo

%
%

\proclaim{\thmlbl{6} Corollary}
Let $E$ be minimizing in $U$ and
$x\in  \partial E$.  Then, in any ball $B(x,r)\subset U$,
the relative volume fractions of $E$ and $E^{c}$ are bounded below by
$\delta(n)>0$.
\endproclaim

\demo{Proof}
By the previous Proposition applied to $E$ and $E^{c}$,
we find that violation of this bound would imply that
$x$ were in the interior of $E$ or of $E^{c}$.  But,
$x\in \partial E$ by assumption, a contradiction.
\enddemo

\proclaim{\thmlbl{7} Corollary }
Let $E$ be minimizing in $U$ and
$x\in  \partial E$.  Then, in any ball $B(x,r)\subset U$,
$P(E,B_r)\geq \delta r^{n-1}$, where $\delta >0$ is an independent
constant.
\endproclaim

\demo{Proof}
This follows from Corollary \thmref{6} plus the explicit form
of the minimizer of $P(A,B_r)$ among sets with
$|A|=|E|$.
\enddemo

{\bf Remark.}
Propositions \thmref{4} and \thmref{5} give an alternative, and more elementary
route to regularity of minimizing sets than the usual path via the
Isoperimetric Theorem for minimal surfaces, cf., \cite{Gi, Chapter 8}.
Using Corollary \thmref{7}, one can go on to show 
existence of tangent cones, etc.   This standard result
is usually proved by reference to the Isometric Theorem for 
minimal surfaces, cf. \cite{Gi, Chapter 5}.

\chapno=7
\subsectionnumber=0
\Increase{\sectionnumber}
\head
\the\chapno. ``Foamy'' sets.
\endhead
\thmno=0
\tagno=0


We conclude by demonstrating existence of sparse, ``foamy''
superminimizing sets having topological boundary with positive
Lebesgue measure,   thus indicating possible limitations of a 
regularity theory for (sub)superminimizing sets.

For $\overline{B}(x_1,r), \overline{B}(x_0,R)\subset U \subset \Bbb{R}^2$,
$\overline{B}(x_1,r)\cap \overline{B}(x_0,R)=\emptyset$,
consider the obstacle problem
$$
\inf \{P(F,U):B(x_1,r)\cup B(x_0,R)\subset F \subset\subset U\}.
\eqnlbl{ob}
$$

\proclaim{\thmlbl{extralemma} Lemma} 
For $r$ sufficiently small, the
solution of \eqnref{ob} is
$$
E=B(x_1,r)\cup B(x_0,R).
$$
Moreover, for any connected
set $\tilde F$ containing
$B(x_1,r)\cup B(x_0,R)$, there holds
$$
P(\tilde F,U)> P(F,U) + \delta,
\eqnlbl{extra}
$$
for some $\delta >0$.
\endproclaim

\demo{Proof}
Without loss of generality, take $U$ to be
all of $\Bbb{R}^2$.  Since we are in two
dimensions, minimal surfaces for \eqnref{ob}
are easily characterized as arcs of $\partial B(x_0,R)$,
$\partial B(x_1,r)$ joined by straight lines.
By explicit comparison, it is then found that
the connected competitor $\tilde F$ with least
perimeter is the convex hull of  $\partial B(x_0,R)$,
$\partial B(x_1,r)$, which for $r$ sufficiently small
satisfies \eqnref{extra}.
Among disconnected competitors, the best is $F=B(x_1,r)\cup B(x_0,R)$,
by \eqnref{iso}.
\enddemo

\proclaim{\thmlbl{foam} Proposition } 
For any open $V\subset \subset U\subset \Bbb{R}^2$, and any
$\varepsilon>0$, there exists a superminimizing set $F$ in $U$ 
such that $\overline{F}=\overline{V}$ and $|F|\le  \pi \varepsilon^2$. 
\endproclaim

\demo{Proof}
Enumerate the rationals as $\{x_j\}$.

{\bf Claim:} For suitably chosen $r_j$, 
$$
F_J:= \cup_{j\le J} B(x_j,r_j)
$$
has the properties:

(i) Any set $F_J\subset G\subset\subset U$
with a connected component containing
two $B(x_j,r_j)$ with $j\le J$, satisfies
$$
P(G,U)> P(F_J,U) + \delta_J, \quad \delta_J>0,
$$

(ii)
$$
\sum_{j=J+1}^{\infty} P(B(x_j,r_j) < \delta_J.
\eqnlbl{deltaii}
$$

{\bf Proof of claim:} The radii $r_j$ may be chosen inductively,
as follows:

Choose $r_1<\varepsilon/2$ sufficiently small that $B(x_1,r_1)\subset V$.
If $x_{j+1}\in F_j$, then take $r_{j+1}=0$.
Otherwise, choose $r_{j+1}$ so small that
$B(x_{j+1},r_{j+1})\subset V\setminus F_j$,
$$
P((B(x_{j+1},r_{j+1}),U)< \delta_j/2,
\eqnlbl{sum}
$$
and, by Lemma \thmref{extralemma}, any connected set
$G$ containing $B(x,r_{j+1}$ and any $B(x_k,r_k)$,
$k\le j$ satisfies
$$
P(G,U)> P(B(x_{j+1},r_{j+1}),U) +
P(B(x_{k},r_{k}),U) + \delta_{j+1}
\eqnlbl{delta}
$$
for some $\delta_{j+1}>0$.  
By \eqnref{sum}, (ii) is clearly satisfied.
Further, \eqnref{sum} and \eqnref{delta} together give (i).
For, if $G$ has a component containing any $B(x_k,r_k)$,
$B(x_l,r_l)$, $k\neq l\le j$, then \eqnref{deltaii} holds
by the induction hypothesis.  Likewise, if no component
of $G$ contains $B(x_{j+1},r_{j+1})$ and any $B(x_k,r_k)$,
$k\le j$.  The remaining case is that precisely one  $B(x_k,r_k)$,
$k\leq j$, lies in a component with $B(x_{j+1},r_{j+1})$, and the rest
lie each in distinct components.  In this case, \eqnref{deltaii}
follows by \eqnref{delta} and \eqnref{iso}.

Defining $F:=\cup_j  B(x_j,r_j)$, we find that
$F$ is superminimizing in $U$.  For, let $G$
be any competitor.  If $G$ has any component
containing $B(x_j,r_j)$ and $B(x_k,r_k)$, $j<k$,
then (i)-(ii) together give 
$$
\eqalign {
P(G,U)&> P(F_k,U) + \delta_k\cr
&> P(F_k,U) + \sum_{k+1}^\infty P(B(x_j,r_j),U) \cr
&\ge P(F,U).
}
$$
On the other hand, if each $B(x_j,r_j)$ lies in a
distinct component $G_j$ of $G$, then either
$G_j\equiv B(x_j,r_j)$, or, by the Isoperimetric
Theorem, $P(G_j,U)\ge P(B(x_j,r_j),U)$, with strict
inequality for some $J$.  Noting that 
$P(G,U)\ge \sum_{j=1}^k P(G_j,U)$ for any finite sum,
and recalling (ii), we thus obtain
$P(G,U)>P(F,U)$ as claimed.

By (ii), and the choice $r_1<\varepsilon$, we have
$|F|\le \pi \varepsilon^2 \sum_{j=1}^\infty (1/2)^{2j} < \pi \varepsilon^2$.
But, clearly, also, $F$ is dense in $V$, giving $\overline{F}=\overline{ V}$
as claimed. 
\qqed
\enddemo

{\bf Remark.}  It is not clear whether such a construction can
be carried out in higher dimensions, since Lemma \thmref{extralemma}
no longer holds with positive $\delta$.

\bigskip

{\bf Consequences:}

1. The construction of Proposition \thmref{foam} shows that in general
$\overline{E}= \overline{E^i}$ is false for subminimizing sets $E$,
in contrast to the result of Corollary \thmref{6}.
It would seem that some form of connectivity must be assumed on $E$, 
if this property is to hold.  

2.  A similar construction with $U=B(0,1)$ yields a superminimizing set
$G$ contained in and dense in  the lower hemisphere $B^-(0,1):=\{x:s\in B(0,1),
\, x_n\leq 0\}$.
Taking $E:=G^{c}$, $F=B^-(0,1)$,  we find
that the strong maximum principle as stated in Theorem \thmref{2.2} is violated.
However, we remark that in the original form as stated in \cite{S2},
the conclusion of the theorem was that $\partial E$ and $\partial F$ should 
agree on their components of $x_0$.  This version of the theorem remains valid
also for the above example, though the two statements are equivalent for
minimizing sets.  Here, again, $\overline{E}\ne \overline{E^i}=\emptyset$,
violating the regularity assumption of \thmref{maxprin}.

Evidently, the issue of a maximum principle for sub- and superminimizing
sets is a delicate one, requiring ideas beyond those in this paper.  This 
would appear to be an interesting area for further study.
\newpage

\head References\endhead

\enddocument

\comment
\ref
\key
\by
\paper
\jour
\vol
\yr
\pages
\endref

\ref
\key
\by
\paper
\jour
\vol
\yr
\pages
\endref
\endcomment

\bye

%% file: inputs.tex
\expandafter \def \csname +pdirichlet\endcsname {1.1}
\expandafter \def \csname -pdirichlet/\endcsname {1}
\expandafter \def \csname +1.2\endcsname {1.2}
\expandafter \def \csname +1.3\endcsname {1.3}
\expandafter \def \csname -1.2/\endcsname {2}
\expandafter \def \csname -1.3/\endcsname {2}
\expandafter \def \csname +2.1\endcsname {2.1}
\expandafter \def \csname +2.2\endcsname {2.2}
\expandafter \def \csname +2.3\endcsname {2.3}
\expandafter \def \csname +2.4\endcsname {2.4}
\expandafter \def \csname +gg\endcsname {2.5}
\expandafter \def \csname -2.1/\endcsname {3}
\expandafter \def \csname -2.2/\endcsname {3}
\expandafter \def \csname -2.3/\endcsname {3}
\expandafter \def \csname +2.6\endcsname {2.6}
\expandafter \def \csname +2.7\endcsname {2.7}
\expandafter \def \csname +2.8\endcsname {2.8}
\expandafter \def \csname +defineE\endcsname {2.9}
\expandafter \def \csname +2.10\endcsname {2.10}
\expandafter \def \csname +coarea\endcsname {2.11}
\expandafter \def \csname -2.4/\endcsname {4}
\expandafter \def \csname -gg/\endcsname {4}
\expandafter \def \csname -2.6/\endcsname {4}
\expandafter \def \csname -2.7/\endcsname {4}
\expandafter \def \csname -2.8/\endcsname {4}
\expandafter \def \csname -defineE/\endcsname {4}
\expandafter \def \csname -2.10/\endcsname {4}
\expandafter \def \csname +reversecoa\endcsname {2.12}
\expandafter \def \csname +iso\endcsname {2.13}
\expandafter \def \csname +2.14\endcsname {2.14}
\expandafter \def \csname +2.15\endcsname {2.15}
\expandafter \def \csname +2.17\endcsname {2.16}
\expandafter \def \csname -coarea/\endcsname {5}
\expandafter \def \csname -reversecoa/\endcsname {5}
\expandafter \def \csname -iso/\endcsname {5}
\expandafter \def \csname -2.14/\endcsname {5}
\expandafter \def \csname -2.15/\endcsname {5}
\expandafter \def \csname +2.18\endcsname {2.17}
\expandafter \def \csname +2.19\endcsname {2.18}
\expandafter \def \csname +2.20\endcsname {2.19}
\expandafter \def \csname +2.21\endcsname {2.20}
\expandafter \def \csname +2.22\endcsname {2.21}
\expandafter \def \csname +2.23\endcsname {2.22}
\expandafter \def \csname +Theorem 2.2\endcsname {3.1}
\expandafter \def \csname +Lemma A\endcsname {3.2}
\expandafter \def \csname -2.17/\endcsname {6}
\expandafter \def \csname -2.18/\endcsname {6}
\expandafter \def \csname -2.19/\endcsname {6}
\expandafter \def \csname -2.20/\endcsname {6}
\expandafter \def \csname -2.21/\endcsname {6}
\expandafter \def \csname -2.22/\endcsname {6}
\expandafter \def \csname -2.23/\endcsname {6}
\expandafter \def \csname +maxprin\endcsname {3.3}
\expandafter \def \csname +1\endcsname {3.1}
\expandafter \def \csname +2\endcsname {3.2}
\expandafter \def \csname +3\endcsname {3.3}
\expandafter \def \csname +4\endcsname {3.4}
\expandafter \def \csname +5\endcsname {3.5}
\expandafter \def \csname +5a\endcsname {3.6}
\expandafter \def \csname -Theorem 2.2/\endcsname {7}
\expandafter \def \csname -Lemma A/\endcsname {7}
\expandafter \def \csname -maxprin/\endcsname {7}
\expandafter \def \csname -1/\endcsname {7}
\expandafter \def \csname -2/\endcsname {7}
\expandafter \def \csname -3/\endcsname {7}
\expandafter \def \csname -4/\endcsname {7}
\expandafter \def \csname +6\endcsname {3.7}
\expandafter \def \csname +7\endcsname {3.8}
\expandafter \def \csname +maxprinsuper\endcsname {3.4}
\expandafter \def \csname -5/\endcsname {8}
\expandafter \def \csname -5a/\endcsname {8}
\expandafter \def \csname -6/\endcsname {8}
\expandafter \def \csname -7/\endcsname {8}
\expandafter \def \csname -maxprinsuper/\endcsname {8}
\expandafter \def \csname +3.1\endcsname {4.1}
\expandafter \def \csname +3.2\endcsname {4.2}
\expandafter \def \csname +lemma3.1\endcsname {4.1}
\expandafter \def \csname +Lemma 3.2\endcsname {4.2}
\expandafter \def \csname -3.1/\endcsname {9}
\expandafter \def \csname -3.2/\endcsname {9}
\expandafter \def \csname -lemma3.1/\endcsname {9}
\expandafter \def \csname +scontain\endcsname {4.3}
\expandafter \def \csname +3.3\endcsname {4.3}
\expandafter \def \csname +2.8a\endcsname {4.4}
\expandafter \def \csname +3.4\endcsname {4.5}
\expandafter \def \csname +3.5\endcsname {4.6}
\expandafter \def \csname -Lemma 3.2/\endcsname {10}
\expandafter \def \csname -scontain/\endcsname {10}
\expandafter \def \csname -3.3/\endcsname {10}
\expandafter \def \csname -2.8a/\endcsname {10}
\expandafter \def \csname -3.4/\endcsname {10}
\expandafter \def \csname +3.7\endcsname {4.7}
\expandafter \def \csname +3.8\endcsname {4.8}
\expandafter \def \csname -3.5/\endcsname {11}
\expandafter \def \csname -3.7/\endcsname {11}
\expandafter \def \csname -3.8/\endcsname {11}
\expandafter \def \csname +3.9\endcsname {4.9}
\expandafter \def \csname +goodcon\endcsname {4.4}
\expandafter \def \csname -3.9/\endcsname {12}
\expandafter \def \csname +3.10\endcsname {4.10}
\expandafter \def \csname +3.11\endcsname {4.11}
\expandafter \def \csname +3.12\endcsname {4.12}
\expandafter \def \csname -goodcon/\endcsname {13}
\expandafter \def \csname -3.10/\endcsname {13}
\expandafter \def \csname -3.11/\endcsname {13}
\expandafter \def \csname -3.12/\endcsname {13}
\expandafter \def \csname +3.17\endcsname {4.13}
\expandafter \def \csname +3.18\endcsname {4.14}
\expandafter \def \csname -3.17/\endcsname {15}
\expandafter \def \csname -3.18/\endcsname {15}
\expandafter \def \csname +gw\endcsname {4.15}
\expandafter \def \csname +lvlset\endcsname {4.16}
\expandafter \def \csname +acontain\endcsname {4.17}
\expandafter \def \csname +defu\endcsname {4.18}
\expandafter \def \csname +dsoln\endcsname {4.5}
\expandafter \def \csname -gw/\endcsname {16}
\expandafter \def \csname -lvlset/\endcsname {16}
\expandafter \def \csname -acontain/\endcsname {16}
\expandafter \def \csname -defu/\endcsname {16}
\expandafter \def \csname -dsoln/\endcsname {16}
\expandafter \def \csname +soln\endcsname {4.6}
\expandafter \def \csname +min\endcsname {4.19}
\expandafter \def \csname +2.22\endcsname {4.20}
\expandafter \def \csname +q\endcsname {4.21}
\expandafter \def \csname +change\endcsname {4.22}
\expandafter \def \csname +star\endcsname {4.23}
\expandafter \def \csname +holderatbdry\endcsname {5.1}
\expandafter \def \csname -soln/\endcsname {17}
\expandafter \def \csname -min/\endcsname {17}
\expandafter \def \csname -2.22/\endcsname {17}
\expandafter \def \csname -q/\endcsname {17}
\expandafter \def \csname -change/\endcsname {17}
\expandafter \def \csname -star/\endcsname {17}
\expandafter \def \csname +lap\endcsname {5.1}
\expandafter \def \csname -holderatbdry/\endcsname {18}
\expandafter \def \csname -lap/\endcsname {18}
\expandafter \def \csname +posmeancurvature\endcsname {5.2}
\expandafter \def \csname +propercontain\endcsname {5.3}
\expandafter \def \csname +goodt\endcsname {5.4}
\expandafter \def \csname +badinequality\endcsname {5.5}
\expandafter \def \csname -posmeancurvature/\endcsname {19}
\expandafter \def \csname -propercontain/\endcsname {19}
\expandafter \def \csname -goodt/\endcsname {19}
\expandafter \def \csname -badinequality/\endcsname {19}
\expandafter \def \csname +holder\endcsname {5.2}
\expandafter \def \csname +cycle\endcsname {5.6}
\expandafter \def \csname -holder/\endcsname {20}
\expandafter \def \csname -cycle/\endcsname {20}
\expandafter \def \csname +psiholder\endcsname {5.7}
\expandafter \def \csname +1\endcsname {6.1}
\expandafter \def \csname +2\endcsname {6.2}
\expandafter \def \csname -psiholder/\endcsname {21}
\expandafter \def \csname -1/\endcsname {21}
\expandafter \def \csname +3\endcsname {6.3}
\expandafter \def \csname +2\endcsname {6.1}
\expandafter \def \csname +3\endcsname {6.2}
\expandafter \def \csname +4\endcsname {6.4}
\expandafter \def \csname -2/\endcsname {22}
\expandafter \def \csname -3/\endcsname {22}
\expandafter \def \csname -2/\endcsname {22}
\expandafter \def \csname -3/\endcsname {22}
\expandafter \def \csname +5\endcsname {6.5}
\expandafter \def \csname +5.1\endcsname {6.6}
\expandafter \def \csname +6\endcsname {6.7}
\expandafter \def \csname +7\endcsname {6.8}
\expandafter \def \csname -4/\endcsname {23}
\expandafter \def \csname -5/\endcsname {23}
\expandafter \def \csname -5.1/\endcsname {23}
\expandafter \def \csname +ob\endcsname {7.1}
\expandafter \def \csname +extralemma\endcsname {7.1}
\expandafter \def \csname +extra\endcsname {7.2}
\expandafter \def \csname +foam\endcsname {7.2}
\expandafter \def \csname -6/\endcsname {24}
\expandafter \def \csname -7/\endcsname {24}
\expandafter \def \csname -ob/\endcsname {24}
\expandafter \def \csname -extralemma/\endcsname {24}
\expandafter \def \csname -extra/\endcsname {24}
\expandafter \def \csname +deltaii\endcsname {7.3}
\expandafter \def \csname +sum\endcsname {7.4}
\expandafter \def \csname +delta\endcsname {7.5}
\expandafter \def \csname -foam/\endcsname {25}
\expandafter \def \csname -deltaii/\endcsname {25}
\expandafter \def \csname -sum/\endcsname {25}
\expandafter \def \csname -delta/\endcsname {25}


\loadbold

\def\upchi{\raise1.2pt\hbox{$\chi$}}		

\def\charfn#1{{\raise1.2pt\hbox{$\chi
_{\kern-1pt\lower3pt\hbox{{$\scriptstyle#1$}}}$}}}

\def\norm#1{\left\Vert{#1}\right\Vert}                 	
\def\abs#1{\left\vert{#1}\right\vert}                   	

\def\today{\ifcase\month\or January\or February\or March\or April\or
May\or June\or July\or August\or September\or October\or November\or
December\fi\space\number\day, \number\year}  

\def\intave#1{\int_{#1}\hbox{\llap{$\raise2.3pt\hbox{\vrule
height.9pt width7pt}\phantom{\scriptstyle{#1}}\mkern-2mu$}}}	

\def\div{{\rm div}}			     	
\def\sqr#1#2{{\vcenter{\hrule height.#2pt\hbox{\vrule
     width.#2pt height#1pt\kern#1pt\vrule width.#2pt}\hrule height.#2pt}}}



 5	
 4		
 3	
 2		
\font\eightpt=cmr8				
 5	
 4	
 3	
 2	

\def\ztextindent#1{\indent\llap{\hbox to\parindent{#1\hfill}}\ignorespaces}
\def\medcup{\mathop{\textstyle\bigcup}\limits}
\def\medcap{\mathop{\textstyle\bigcap}\limits}
\def\nullset{\hbox{\eightpt \O}}

\def\b+{B^{+}(x_{0},R)}

\def\sumh_#1^#2{\sum\limits_{#1}^{#2}}

\def\-{\bar}
\def\~{\tilde}
\def\){\bigl)}
\def\({\bigl(}

\predefine\uu{\u}
\def\u{\bar{u}}
\predefine\vv{\v}
\def\v{\bar{v}}

\def\bar#1{\overline{#1}}

\def\qqed{\enspace\hfill $\square$}

\def\osc{\operatorname{osc}}
\def\doscball#1#2{\mathop{\osc\;#1}_{\kern-13pt B(x_{0},#2)}}

\def\R{\bold R}

\def\rn{\bold R^n}
\def\r1{\bold R^1}

\redefine\div{\operatorname{div}}

\define\fn#1\endfn{$\bigstar$\footnote {#1} }


\define\Cpr#1{\boldsymbol \gamma_{p;#1}}
\define\Cprn#1{\boldsymbol\gamma_{n;#1}}

\def\hidebr#1{} 

\medskip 

\newcount\sectionnumber
	\sectionnumber=0
	
\newcount\subsectionnumber
	\subsectionnumber=0



\font\sixrm=cmr6
\newcount\tagno \tagno=0		        
\newcount\thmno	\thmno=0	         	
\newcount\bibno	\bibno=0			
\newcount\chapno\chapno=0                       
\newif\ifproofmode
\proofmodetrue
\newif\ifwanted
\wantedfalse
\newif\ifindexed
\indexedfalse

\def\ifundefined#1{\expandafter\ifx\csname+#1\endcsname\relax}

\def\Wanted#1{\ifundefined{#1} \wantedtrue
\immediate\write0{Wanted #1\the\chapno.\the\thmno}\fi}

\def\Increase#1{{\global\advance#1 by 1}}

\def\Assign#1#2{\immediate\write1{\noexpand\expandafter\noexpand\def
 \noexpand\csname+#1\endcsname{#2}}\relax
 \global\expandafter\edef\csname+#1\endcsname{#2}}

\def\pAssign#1#2{\write1{\noexpand\expandafter\noexpand\def
 \noexpand\csname-#1\endcsname{#2}}}

\def\lPut#1{\ifproofmode\llap{\hbox{\sixrm #1\ \ \ }}\fi}
\def\rPut#1{\ifproofmode$^{\hbox{\sixrm #1}}$\fi}



\def\chp#1{\global\tagno=0\global\thmno=0\Increase\chapno\Assign{#1}
{\the\chapno}{\lPut{#1}\the\chapno}}


\def\thm#1{\Increase\thmno
  \Assign{#1}{\the\chapno.\the\thmno}
  \lPut{#1}{\the\chapno.\the\thmno}}


\def\frm#1{\Increase\tagno
  \Assign{#1}{\the\chapno.\the\tagno}
  \lPut{#1}{\the\chapno.\the\tagno}}


\def\bib#1{\Increase\bibno\Assign{#1}{\the\bibno}
\lPut{#1}{\the\bibno}}


\def\pgp#1{\pAssign{#1/}{\the\pageno}}


\def\ix#1#2#3{\pAssign{#2}{\the\pageno}
\immediate\write#1{\noexpand\idxitem{#3}{\noexpand\csname+#2\endcsname}}}


\def\rf#1{\Wanted{#1}\csname+#1\endcsname\relax\rPut{#1}}


\def\rfp#1{\Wanted{#1}\csname-#1/\endcsname\relax\rPut{#1}}

\immediate\openout1=\jobname.aux

\ifindexed
\immediate\openout2=\jobname.idx
\immediate\openout3=\jobname.sym
\fi


\def\thmref#1{{\rm{\rf{#1}}}}
\def\eqnref#1{{\rm{(\rf{#1})}}}
\def\markp#1{{\rm{\pgp{#1}}}}
\def\referp#1{{\kern-3pt\rm{\rfp{#1}}}}
\def\thmlbl#1{\thm{#1}\markp{#1}}
\def\eqnlbl#1{\tag{\kern-3pt\frm{#1}}\markp{#1}}

\input amsppt.sty
\catcode`@=11


\def\roster{%
  \envir@stack\endroster
 \edef\leftskip@{\leftskip\the\leftskip}%
 \relaxnext@
 \rostercount@\z@
 \def\romitem{\FN@\rosteritem@}
 \DN@{\ifx\next\runinitem\let\next@\nextii@\else
  \let\next@\nextiii@\fi\next@}%
 \DNii@\runinitem
  {\unskip
   \DN@{\ifx\next[\let\next@\nextii@\else
    \ifx\next"\let\next@\nextiii@\else\let\next@\nextiv@\fi\fi\next@}%
   \DNii@[####1]{\rostercount@####1\relax
    \enspace\therosteritem{\romannumeral\rostercount@}~\ignorespaces}%
   \def\nextiii@"####1"{\enspace{\rm####1}~\ignorespaces}%
   \def\nextiv@{\enspace\therosteritem1\rostercount@\@ne~}%
   \par@\firstitem@false
   \FN@\next@}
 \def\nextiii@{\par\par@
  \penalty\@m\smallskip\vskip-\parskip
  \firstitem@true}
 \FN@\next@}
\def\rosteritem@{\iffirstitem@\firstitem@false
  \else\par\vskip-\parskip\fi
 \leftskip\rosteritemwd \advance\leftskip\normalparindent
 \advance\leftskip.5em \noindent
 \DNii@[##1]{\rostercount@##1\relax\itembox@}%
 \def\nextiii@"##1"{\def\therosteritem@{\rm##1}\itembox@}%
 \def\nextiv@{\advance\rostercount@\@ne\itembox@}%
 \def\therosteritem@{\therosteritem{\romannumeral\rostercount@}}%
 \ifx\next[\let\next@\nextii@\else\ifx\next"\let\next@\nextiii@\else
  \let\next@\nextiv@\fi\fi\next@}
\def\Runinitem#1\roster\runinitem{\relaxnext@
  \envir@stack\endroster
 \rostercount@\z@
 \def\romitem{\FN@\rosteritem@}%
 \def\runinitem@{#1}%
 \DN@{\ifx\next[\let\next\nextii@\else\ifx\next"\let\next\nextiii@
  \else\let\next\nextiv@\fi\fi\next}%
 \DNii@[##1]{\rostercount@##1\relax
  \def\romitem@{\therosteritem{\romannumeral\rostercount@}}\nextv@}%
 \def\nextiii@"##1"{\def\romitem@{{\rm##1}}\nextv@}%
 \def\nextiv@{\advance\rostercount@\@ne
  \def\romitem@{\therosteritem{\romannumeral\rostercount@}}\nextv@}%
 \def\nextv@{\setbox\z@\vbox
  {\ifnextRunin@\noindent\fi
  \runinitem@\unskip\enspace\romitem@~\par
  \global\rosterhangafter@\prevgraf}%
  \firstitem@false
  \ifnextRunin@\else\par\fi
  \hangafter\rosterhangafter@\hangindent3\normalparindent
  \ifnextRunin@\noindent\fi
  \runinitem@\unskip\enspace
  \romitem@~\ifnextRunin@\else\par@\fi
  \nextRunin@true\ignorespaces}
 \FN@\next@}